\RequirePackage{fix-cm}
\documentclass[smallextended,envcountsect,envcountsame]{svjour3}
\smartqed
\usepackage[T1]{fontenc}
\usepackage{lmodern}
\usepackage{amsmath,amssymb,mathtools}
\usepackage{graphicx,booktabs,array,tabularx}
\usepackage{microtype}
\usepackage[numbers,sort&compress]{natbib}
\usepackage[hidelinks]{hyperref}
\hypersetup{pdftitle={},pdfauthor={},pdfsubject={},pdfkeywords={},
 pdfcreator={},pdfproducer={}}
\ifdefined\pdfinfoomitdate\pdfinfoomitdate=1\fi
\ifdefined\pdftrailerid\pdftrailerid{}\fi
\ifdefined\pdfsuppressptexinfo\pdfsuppressptexinfo=15\fi
\makeatletter
\g@addto@macro\UrlBreaks{\do\a\do\b\do\c\do\d\do\e\do\f
 \do0\do1\do2\do3\do4\do5\do6\do7\do8\do9}
\makeatother
\providecommand{\doi}[1]{}
\renewcommand{\doi}[1]{\url{https://doi.org/#1}}
\numberwithin{equation}{section}
\newcommand{\F}{\mathbb F}
\newcommand{\R}{\mathbb R}
\newcommand{\C}{\mathbb C}
\newcommand{\E}{\mathbb E}
\newcommand{\T}{\mathcal T}
\DeclareMathOperator{\tr}{tr}
\DeclareMathOperator{\rank}{rank}
\DeclareMathOperator{\OPT}{OPT}
\DeclareMathOperator{\St}{St}
\DeclareMathOperator{\diag}{diag}
\DeclareMathOperator{\Mean}{Mean}
\newcommand{\factor}[3][]{\alpha^{#1}_{#2}(#3)}
\newcommand{\directfactor}[3]{\alpha_{#1,#2}(#3)}
\newcommand{\eps}{\varepsilon}
\newcommand{\dd}{\,\mathrm d}
\journalname{}
\renewcommand{\makeheadbox}{}
\begin{document}
\title{Sharp approximation thresholds for power means of common-subspace volumes}
\titlerunning{Sharp thresholds for common-subspace volumes}
\author{Alexey Peregudin}
\authorrunning{Alexey Peregudin}
\institute{Alexey Peregudin \at
School of Electrical and Electronic Engineering, University of Sheffield\\
\email{a.peregudin@sheffield.ac.uk}}
\date{}
\maketitle
\begin{abstract}
We study how well a single $k$-dimensional subspace can preserve volume
from several positive semidefinite matrices at once. Each matrix is
compressed to the subspace and scored by the geometric mean of the
eigenvalues of that compression; the scores are aggregated by a weighted
power mean of order $p \le 1$. For every fixed $k$ and every fixed rational
$p$ we determine the optimal approximation factor: it is achieved by a
randomized polynomial-time algorithm, and every larger constant is
NP-hard. Above a critical power ($p=-k/2$ over the reals and $p=-k$ over
the complex numbers) the factor is an explicit product of Gamma functions.
At and below the critical power, every positive constant-factor
approximation is NP-hard. The factor comes from a concave relaxation over
the positive semidefinite matrices of trace one, which agrees with the
volume score at every subspace, followed by ordinary Gaussian rounding.
For two rank-$k$ matrices the relaxation is exact, giving a closed-form
optimum for Martin means of two subspaces.
\keywords{Approximation algorithms \and Determinant maximization \and
Semidefinite relaxation \and Gaussian rounding \and
Hardness of approximation \and Grassmannian optimization}
\subclass{90C26 \and 90C22 \and 68Q17 \and 15B52}
\end{abstract}

\section{Introduction}\label{sec:introduction}

A common subspace can retain substantial energy from several observations
while losing an entire direction in one of them. A determinant detects
this loss: it vanishes as soon as the compressed observation becomes
singular. This makes volume preservation a natural objective when every
direction of a selected subspace must remain represented.

Let $A_1,\ldots,A_m$ be positive semidefinite matrices and let the columns
of $P$ form an orthonormal basis for a $k$-dimensional subspace. We seek
to maximize
\[
 \Mean_{p,w}\bigl[\det(P^*A_iP)^{1/k}\bigr],
 \qquad P^*P=I_k,
\]
where $\Mean_{p,w}$ is the weighted power mean across observations.
The $k$th root makes the score homogeneous of degree one in $A_i$,
as for a quadratic form $x^*A_ix$. Arithmetic averaging rewards the average score;
geometric and negative-power means penalize poorly represented
observations more strongly.

One motivating example is the Martin mean of subspaces
\cite{YangHospedales2016}. For observed orthonormal frames $U_i$, its
loss is the weighted sum of $-\log|\det(U_i^*P)|^2$. Each absolute determinant
is the product of the cosines of the principal angles, so an orthogonal
direction incurs infinite loss. This is our geometric-mean problem
with $A_i\coloneqq U_iU_i^*$. Local manifold methods are available; the question
here is what a polynomial-time algorithm can guarantee globally.

We give two answers to the question of how much is lost when a
solution of a matrix relaxation is rounded to one subspace.
For general instances, we determine the optimal approximation factor
independent of the ambient dimension, for every fixed selected
dimension and rational power at most one, conditional on
$\mathrm{NP}\not\subseteq\mathrm{BPP}$.
Above the critical power $-k/2$ over $\R$ and $-k$ over $\C$,
an explicit Gaussian determinant moment gives the factor.
At and below this power, every positive
constant-factor approximation is NP-hard.
The algorithm attains the positive factor, with a qualification at
$(k,p)=(1,1)$ imposed by rational output.
Matching hardness at $p=1$ requires a tensor construction
(Section~\ref{sec:hardness}).

For two rank-$k$ observations, the answer is stronger: there is no
rounding loss. Every positive semidefinite matrix of trace $k$ can be
replaced by one subspace without decreasing either determinant score.
This two-view purification holds before choosing an aggregate,
so it gives exact reformulations for all coordinatewise increasing
aggregates. For power means of order at most one, the resulting
concave formulation also gives efficient relative approximation.

Both results use the same determinant-sensitive trace relaxation.
Its domain is $\T_n\coloneqq \{X\succeq0:\tr X=1\}$. The relaxed score
of an observation $A$ is
\[
 f_A(X)\coloneqq k\,e_k(A^{1/2}XA^{1/2})^{1/k},
 \qquad f_A(PP^*/k)=v_A(P),
\]
where $e_k$ is the elementary symmetric polynomial of the eigenvalues.
This score is concave. When $A=UDU^*$ has rank $k$, with
$U^*U=I_k$ and $d\coloneqq \det(D)^{1/k}$, setting $Y\coloneqq kX$ gives
$f_A(X)=d\det(U^*YU)^{1/k}$, precisely the quantity preserved by
the two-view construction. Independent Gaussian columns give the
general approximation; principal-angle blocks give lossless rounding
in the two-view case.

Martin means illustrate both conclusions. The geometric-power
approximation becomes a sharp dimension-free additive bound for the
Martin loss. For two observed subspaces, purification
reduces the problem to scalar maximizations in their principal-angle
blocks, giving an explicit global solution.

\subsection{Relation to earlier work}\label{sec:related}

\paragraph{Known scalar factors and Yuan's subspace proposal.}
Barvinok \cite{Barvinok2014} studied approximation of quadratic images
by their convex hulls in relative entropy.
Yuan and Parrilo's 2021 preprint \cite{YuanParrilo2021Products}
gives the rank-sensitive scalar geometric-mean factor
(Theorem~1.3) and matching relaxation gaps (Theorem~5.1).
These are also Theorems~4.1.3 and~4.5.1 of Yuan's dissertation
\cite{YuanThesis2022}.
Its Theorem~5.4.2 gives the scalar positive-power factors for
$0<p\le1$: after taking the $1/p$ power, they are exactly our
$\directfactor{r}{1}{p}$.
The scalar factors and Gaussian covariance rounding are therefore
established predecessors of the present results.

More directly, Section~5.6 of the dissertation formulates the
common-subspace log-determinant problem, derives a trace relaxation,
and proposes independent Gaussian columns followed by polar
normalization. It leaves the approximation analysis for future work.
That trace relaxation can lose an arbitrarily large factor for
determinant volume, even for one positive-definite observation
(Section~\ref{model:trace-gap}). Our concave relaxed score instead
agrees with the determinant score at every subspace.
We establish its sharp approximation and hardness thresholds across
the power-mean scale.

The related PSD-permanent literature includes the covariance
relaxation of Anari et al.\ \cite{AnariEtAl2017Permanent} and
the product-of-linear-forms results of Yuan and Parrilo
\cite{YuanParrilo2022Permanent}.
Anari and Ebrahimnejad's 2026 preprint
\cite{AnariEbrahimnejad2026Permanent} obtains the optimal exponential approximation
exponent, assuming $\mathrm P\ne\mathrm{NP}$, through an
entropy-corrected concave covariance relaxation.
Its objective is a permanent, rather than a determinant of a
compression to a common subspace.

\paragraph{Gaussian rounding and sharp constants.}
Gaussian rounding is central to semidefinite approximation
\cite{Nesterov1998,SoZhangYe2007}.
For matrix objectives, Nemirovski \cite{Nemirovski2007} and
Naor, Regev and Vidick \cite{NaorRegevVidick2014} developed
semidefinite approximation under orthogonality constraints and
noncommutative Grothendieck rounding.
Bandeira, Kennedy and Singer \cite{BandeiraKennedySinger2016}
combine Gaussian matrix rounding and polar projection in the real and
complex settings. Their constants involve mean Gaussian singular
values, with matching relaxation gaps. Their objective is quadratic
in several orthogonal or unitary blocks; our determinant-root means
require determinant moments.
Bri\"et, de Oliveira Filho and Vallentin
\cite{BrietOliveiraVallentin2010} obtain Gamma-ratio approximation
constants for the PSD Grothendieck problem with a Gram-rank
constraint, with optimality under the Unique Games Conjecture.
That Gram rank differs from the covariance rank in our rounding bound.

Sharp NP-hardness thresholds governed by Gaussian moments also occur
in geometric approximation \cite{GuruswamiEtAl2016}.
Our reductions use the real and complex norm-hardness results
of Bri\"et, Regev and Saket \cite{BrietRegevSaket2017} and the
block-subspace promise in the full version of
Ebrahimnejad, Nagda and Oveis Gharan
\cite{EbrahimnejadNagdaOveisGharan2024}.
The Yuan--Parrilo product-of-forms hardness theorem
\cite[Theorem~7.1]{YuanParrilo2021Products} has factor
$(1-\varepsilon)^{1/d}$ for $d=\Omega(n)$ forms
(also \cite[Theorem~4.7.1]{YuanThesis2022});
it does not identify a sharp constant-factor threshold.
Tensor amplification also appears in norm hardness
\cite{BhattiproluEtAl2023}; our arithmetic-power argument controls
coefficient energies across tensor orientations.
General Stiefel and Grassmann hardness
\cite{LaiLimTang2026,LaiLimYe2025} concerns broader objective classes.

\paragraph{Elementary-symmetric volume criteria.}
The concavity of elementary-symmetric roots
\cite{MarcusLopes1957,Garding1959,BauschkeEtAl2001} and
Gaussian determinant moments
\cite{Muirhead1982,Goodman1963Determinant,EdelmanRao2005,Nielsen1999}
are established ingredients.
Pronzato, Wynn and Zhigljavsky \cite{PronzatoWynnZhigljavsky2017}
interpret elementary symmetric polynomials of covariance eigenvalues
through expected squared simplex volumes and study their concave roots.
Mariet and Sra \cite{MarietSra2017} use elementary symmetric
polynomials as partial-volume criteria in experimental design.
These give geometric and statistical context for our relaxed score
$e_k^{1/k}$. Hyperbolic optimization supplies a broader computational
framework \cite{Guler1997}; Saunderson and Parrilo
\cite{SaundersonParrilo2015} give polynomial-size semidefinite
representations of the associated derivative cones.
Our algorithm uses the evaluation and separation oracles proved in
Appendix~\ref{app:algorithms}; a cone representation alone does not
formulate the full power-mean objective.

\paragraph{Two-observation exactness.}
Fair dimensionality reduction
\cite{SamadiEtAl2018,TantipongpipatEtAl2019}
studies several retained-energy criteria.
Tantipongpipat et al.\ \cite[Theorem~1.2]{TantipongpipatEtAl2019}
prove exactness for two affine criteria under a concave aggregation
nondecreasing in at least one argument.
Li and Xie \cite[Corollary~4]{LiXie2026} recover two-criterion
fair-PCA exactness within a general theory of partial convexification
and rank bounds. Their fair-PCA relaxation includes an eigenvalue cap.
Shen et al.\ \cite{ShenDavisLuBai2026} exploit hidden convexity
to solve two-group fair PCA through eigenvalue optimization of
reconstruction-loss trace criteria.
Randomized SDP rank reduction \cite{SoYeZhang2008} controls relative
deviations in linear measurements.

Our two-view purification theorem instead simultaneously preserves
two nonlinear determinant scores from the full trace-$k$ PSD domain,
under the assumption that both observations have rank $k$.
Its principal-angle blocks come from classical two-subspace geometry
\cite{Halmos1969,BjorckGolub1973}.
The author's 2026 preprint \cite{Peregudin2026Majorants}
gives scalar exactness under a lower-envelope convexity condition,
automatic for two real forms (Proposition~3.2 and
Corollaries~4.1--4.2).
Those scalar results allow general observation ranks; the
two-view determinant theorem here assumes rank-$k$ observations.

\paragraph{Subspace means and optimization methods.}
Stiefel and Grassmann optimization has a classical geometric
foundation \cite{EdelmanAriasSmith1998}; a recent account is
\cite{BendokatZimmermannAbsil2024}.
Feasible manifold methods \cite{WenYin2013} and nonlinear
eigenvalue methods \cite{Li2026NEPv} provide iterative algorithms;
the convergence guarantees in the latter framework are to stationary
points under its assumptions.
Semidefinite methods for quadratic frame objectives
\cite{BurerPark2024,GilmanBurerBalzano2025,CoryWrightPauphilet2026}
also relax orthogonality. Cory-Wright and Pauphilet's extended
version \cite[arXiv version~4]{CoryWrightPauphilet2026}
gives a $1/3$ approximation using Gaussian polar rounding.
These objectives are quadratic in a frame; ours are generally
nonlinear and invariant under changes of basis.
Xu's 2026 preprint \cite{Xu2026Fair}
uses a smooth minimum of compressed log determinants, equivalently
a logarithmic negative-power-mean objective, and develops
stationary-point methods for the continuous subspace stage.

The Martin loss originates in a metric for ARMA processes
\cite{Martin2000}, with its subspace-angle interpretation developed
in \cite{DeCockDeMoor2002}.
Yang and Hospedales \cite{YangHospedales2016} explicitly optimize
the finite-dimensional Martin mean. Our results give a global
approximation guarantee and an explicit solution for two subspaces.
The relaxed weighted log-determinant problem belongs to the MAXDET
and weighted determinant-maximization families
\cite{VandenbergheBoydWu1998,TsuchiyaXia2007}.
Binet--Cauchy kernels \cite{VishwanathanSmolaVidal2007} use
determinant overlaps; averaging squared overlaps corresponds to
power $p=k$ in our scale. Intrinsic Grassmann means
\cite{BatziesEtAl2015} instead minimize squared geodesic distances,
which gives a different loss from Martin's logarithmic determinant.

\paragraph{Discrete determinant selection.}
Subdeterminant maximization and largest-simplex approximation
\cite{Nikolov2015,NikolovSingh2016,MadanEtAl2020,EbrahimiStraszakVishnoi2017}
select vectors from a finite collection.
Brown et al.\ \cite{BrownEtAl2022} develop constructive
matroid-intersection algorithms for determinant maximization.
Experimental design \cite{SinghXie2020,NikolovSinghTantipongpipat2019}
aggregates selected information matrices; the partial-volume criteria
of \cite{MarietSra2017} connect this literature to elementary
symmetric polynomials.
Here every fixed observation is compressed to the same variable
subspace. Constants depend on normalization: a factor $b$ for a
determinant root is $b^k$ for the determinant.

\paragraph{Organization.}
Section~\ref{sec:formulation} fixes notation and the computational model.
Section~\ref{sec:guarantees} states the classification and its
uniform-input algorithmic form. Section~\ref{sec:lifts} develops the
trace relaxation and Gaussian rounding; Section~\ref{sec:hardness}
gives matching hardness and its transfer to subspace observations.
Section~\ref{sec:observations} proves two-view
purification, and Section~\ref{sec:martin} derives the Martin
consequences. Full technical proofs and numerical guarantees are in
Appendices~\ref{app:rounding}--\ref{proof:exactness}.

\section{Notation and problem formulation}\label{sec:formulation}

\subsection{Notation}
The field $\F$ is $\R$ or $\C$. The symbol $A^*$ denotes transpose
or conjugate transpose, respectively. All observation matrices are
self-adjoint. We write $A\succeq0$ for positive semidefiniteness,
$A\succ0$ for positive definiteness, and $A\succeq B$ when $A-B\succeq0$.
The square root $A^{1/2}$ is the positive semidefinite square root;
it exists also when $A$ is singular. Vector norms are Euclidean,
$\|\cdot\|_{\rm op}$ is the matrix operator norm, and
$\|\cdot\|_{\rm F}$ is the Frobenius norm.

An orthonormal $k$-frame is a matrix in
$\St(n,k)\coloneqq \{P\in\F^{n\times k}:P^*P=I_k\}$.
Its columns span a $k$-dimensional subspace and $PP^*$ is the
orthogonal projector onto that subspace. We write $P_\perp$ for
an orthonormal frame of the orthogonal complement.
The collection of these subspaces is the Grassmannian.

For a positive semidefinite $d\times d$ matrix $C$ with eigenvalues
$\lambda_1,\ldots,\lambda_d$, define
\[
 e_j(C)\coloneqq \sum_{1\le i_1<\cdots<i_j\le d}
          \lambda_{i_1}\cdots\lambda_{i_j},\qquad e_0(C)\coloneqq 1.
\]
Thus $e_1(C)=\tr C$, $e_d(C)=\det C$, and
$\det(I_d+zC)=\sum_{j=0}^d e_j(C)z^j$.
Equivalently, $e_k(C)$ is the sum of the $k\times k$ principal minors of $C$.

We use natural logarithms and the functions
\[
 \Gamma(x)\coloneqq \int_0^\infty t^{x-1}e^{-t}\,\dd t,\qquad
 \psi(x)\coloneqq \frac{d}{dx}\log\Gamma(x)\quad(x>0).
\]
Set $\beta\coloneqq 1/2$ over $\R$ and $\beta\coloneqq 1$ over $\C$.
A standard real Gaussian has variance one; a standard complex
Gaussian is $(\xi+i\zeta)/\sqrt2$ for independent real standard
Gaussians $\xi,\zeta$. A standard Gaussian matrix has independent
such entries, each with expected squared modulus one.
Expectation and probability are denoted by $\E$ and $\Pr$.

\subsection{Volume scores and their means}
Fix $1\le k\le n$, $n\times n$ matrices $A_1,\ldots,A_m\succeq0$,
and positive weights $w_i$ summing to one.
We call the $A_i$ \emph{observations}; in applications they may be
covariance or second-moment matrices, or orthogonal projections onto
observed subspaces.
For nonnegative scalars, the weighted power mean is
\begin{equation}\label{intro:means}
 \Mean_{p,w}[a_1,\ldots,a_m]\coloneqq 
 \begin{cases}
  (\sum_{i=1}^m w_i a_i^p)^{1/p},&p\ne0,\\
  \prod_{i=1}^m a_i^{w_i},&p=0.
 \end{cases}
\end{equation}
At zero coordinates we use the continuous extension: if any
positive-weight coordinate vanishes, the mean is zero for $p\le0$.
Every power in this paper is finite. The notation
$\Mean_{p,w}[a_i]$ abbreviates the list over the displayed observation
index; for a random variable, $\Mean_p$ denotes the corresponding
expectation-based mean.

The compression $P^*AP$ records an observation on the selected
subspace. Its score is the geometric mean of its $k$ eigenvalues:
\begin{equation}\label{intro:problem}
 \begin{aligned}
 v_A(P)&\coloneqq \det(P^*AP)^{1/k}
       =\left(\prod_{j=1}^k\lambda_j(P^*AP)\right)^{1/k},\\
 \OPT_{p,k}(A,w)&\coloneqq \max_{P\in\St(n,k)}\Mean_{p,w}[v_{A_i}(P)].
 \end{aligned}
\end{equation}
Replacing $P$ by $PU$, for an orthogonal or unitary $U$, leaves every
score unchanged. Over $\R$, $\det(P^*AP)$ is the squared
$k$-dimensional volume multiplier under $A^{1/2}$.
The score is homogeneous of degree one in $A$.

There are two distinct averages: the geometric mean within one
compressed observation and the chosen power mean across observations.
The identity
\[
 \Mean_{p,w}[d_i^{1/k}]^k=\Mean_{p/k,w}[d_i]
\]
also shows why arithmetic averaging of raw determinants corresponds
to $p=k$, rather than $p=1$.

\subsection{Input, output, and guarantees}
Computational statements use rational matrix entries, rational powers,
and positive rational weights. Over $\C$, both real and imaginary
parts are rational. A subspace is returned as the range of a
full-column-rank rational matrix $B_{\rm out}$; its rational projector is
$B_{\rm out}(B_{\rm out}^*B_{\rm out})^{-1}B_{\rm out}^*$.
An orthonormal basis for that same subspace need not be rational.

A factor $\rho\in[0,1]$ means a feasible value at least
$\rho\OPT_{p,k}$. Randomized success means probability at least $2/3$;
independent repetition amplifies it. Polynomial time refers to the
complete input bit length, with any further accuracy dependence stated
explicitly.
Hardness fixes the field, selected dimension $k$, rational power $p$,
and target approximation factor $\rho$ while the ambient dimension and
input size grow.
The conditional optimality of randomized algorithms uses
$\mathrm{NP}\not\subseteq\mathrm{BPP}$; the gap reductions are
unconditional NP-hardness statements.

\section{Sharp approximation guarantees}\label{sec:guarantees}

\subsection{The computational threshold}
For the fixed field parameter $\beta$, define, for finite $p\le1$,
\begin{equation}\label{intro:binfty}
 \factor{k}{p}\coloneqq 
 \begin{cases}
 0,&p\le-k\beta,\\[1mm]
 \displaystyle\frac1{k\beta}
 \left[\prod_{j=1}^k\frac{\Gamma(j\beta+p/k)}{\Gamma(j\beta)}
 \right]^{1/p},&p>-k\beta,\quad p\ne0,\\[3mm]
 \displaystyle\frac1{k\beta}
 \exp\left(\frac1k\sum_{j=1}^k\psi(j\beta)\right),&p=0.
 \end{cases}
\end{equation}
For $p>-k\beta$, this is the power mean of
$|\det H|^{2/k}/k$ for a square standard Gaussian $H$.
The value at zero is its geometric mean.
The zero extension at $p=-k\beta$ reflects the computational theorem
below; it is not an evaluation of a divergent Gamma product.

\begin{theorem}[Sharp approximation classification]\label{intro:classification}
Fix the field, an integer $k\ge1$, and a rational power $p\le1$.
For problem~\eqref{intro:problem}:
\begin{enumerate}
\item If $p>-k\beta$ and $(k,p)\ne(1,1)$, a randomized polynomial-time
algorithm attains value at least $\factor{k}{p}\OPT_{p,k}$.
\item For every fixed $\rho$ with $\factor{k}{p}<\rho<1$,
obtaining a $\rho$-approximation is NP-hard.
If $p\le-k\beta$, every positive constant-factor approximation is NP-hard.
\end{enumerate}
The hardness holds for equally weighted rational PSD observations
of rank at most $k$, and also for equally weighted positive-definite
observations as a separate family.
It transfers to rational rank-$k$ projection observations with
positive rational weights; at $p=0$, equal weights suffice.
At $(k,p)=(1,1)$ the problem is a largest-eigenvalue problem,
with approximation supremum one.
\end{theorem}
The proof is assembled in Section~\ref{app:attainment}.
Under $\mathrm{NP}\not\subseteq\mathrm{BPP}$, the theorem determines
the optimal factor independent of $n$ and $m$.
We call $p=-k\beta$ the \emph{negative-moment boundary}, since it is
where the inverse-power Gaussian determinant moment ceases to be finite.
A zero universal factor does not say that all instance optima vanish:
it excludes one fixed positive relative guarantee over the entire family.

\subsection{A single algorithm with dimension and power as inputs}
Theorem~\ref{intro:classification} fixes $k$ and $p$ in advance.
Its algorithmic assertion strengthens to a single algorithm taking both
as inputs. Fix a relative margin $0<\delta<1$ above the negative-moment boundary:
\begin{equation}\label{intro:power-range}
 -(1-\delta)k\beta\le p\le1.
\end{equation}
\begin{theorem}[Attainment of the threshold]\label{alg:attainment}
Fix the field and $0<\delta<1$. One randomized algorithm, polynomial
in the complete rational input bit length, attains
$\factor{k}{p}$ for all input $k,p$ in \eqref{intro:power-range},
except at $(k,p)=(1,1)$.
At $(k,p)=(1,1)$, every prescribed factor
$1-\eps$ is available in time polynomial also in
$\max(1,1/\eps)$; an exact optimum has an algebraic representation.
\end{theorem}
The proof is in Section~\ref{app:attainment}, with the scalar
near-one estimate in Section~\ref{app:scalar-attainment}.
Attainment means that the displayed factor has no additional
subtracted error allowance. The Gamma constant appears only in the
analysis; rational-output limitations at $(k,p)=(1,1)$ are discussed
in Section~\ref{app:attainment}.
For rational $p<1$ arbitrarily close to one, the algorithm remains
polynomial in the bit length of $p$: its sample count has no
factor proportional to $1/(1-p)$.

Table~\ref{intro:quantifiers} distinguishes this assertion from
the finite-dimensional guarantee of Section~\ref{sec:lifts} and the
hardness theorem.

\begin{table}[t]
\caption{The parameters held fixed in the approximation statements.}
\label{intro:quantifiers}
\small
\begin{tabularx}{\textwidth}{@{}>{\raggedright\arraybackslash}p{.24\textwidth}
>{\raggedright\arraybackslash}p{.25\textwidth}
>{\raggedright\arraybackslash}X@{}}
\toprule
Statement & Fixed parameters & Inputs and running-time dependence\\
\midrule
Attained computational factor
& Field, margin $\delta$
& $n,m,k,p,A,w$; polynomial in full bit length.\\
Finite-dimensional factor with slack
& Field, margin $\delta$
& The same data and $\eps>0$; polynomial also in $1/\eps$.\\
Hardness above the threshold
& Field, selected $k$, rational $p$, target factor $\rho$
& Ambient $n$, observation count, and rational data grow.\\
\bottomrule
\end{tabularx}
\end{table}

\section{The determinant relaxation and Gaussian rounding}
\label{sec:lifts}\label{relax:lifts}
A subspace determines a normalized projection matrix. We relax that
matrix to an arbitrary PSD matrix of trace one, while retaining the
determinant information in every observation. Gaussian sampling then
converts the fractional solution back to one subspace.

\subsection{A determinant-sensitive trace relaxation}\label{model:trace-gap}
Let
\begin{equation}\label{relax:domain}
 \T_n\coloneqq \{X\succeq0:\tr X=1\}.
\end{equation}
For $A\succeq0$, define
\begin{equation}\label{relax:direct-lift}
 f_A(X)\coloneqq k\,e_k(A^{1/2}XA^{1/2})^{1/k}.
\end{equation}
We call $f_A$ the \emph{relaxed score} of observation $A$.
The elementary-symmetric-root inequality of Marcus and Lopes
\cite{MarcusLopes1957} makes this function concave on the PSD cone.
It is nonnegative, Loewner-increasing, and homogeneous of degree one.
Monotonicity follows from eigenvalue monotonicity. These properties
extend to the boundary by continuity. If $\rank A<k$, then $f_A$
vanishes identically.

For an orthonormal $k$-frame $P$, the nonzero eigenvalues of
$A^{1/2}PP^*A^{1/2}$ are those of $P^*AP$. Consequently
\begin{equation}\label{relax:direct-exact}
 f_A(PP^*/k)=\det(P^*AP)^{1/k}=v_A(P).
\end{equation}
Since the power mean is concave and coordinatewise increasing for
$p\le1$, the function
\[
 F(X)\coloneqq \Mean_{p,w}[f_{A_i}(X)]
\]
is concave on $\T_n$, with the stated continuous convention at zero
scores. Maximizing it is therefore a convex optimization problem,
and its optimum $S^\star$ satisfies $S^\star\ge\OPT_{p,k}$.
The convex hull of normalized rank-$k$ projectors additionally imposes
$X\preceq I/k$ \cite{OvertonWomersley1992}; the domain
\eqref{relax:domain} omits this spectral cap.

\paragraph{Why the relaxed score must retain determinants.}
The trace constraint does not make this a linear trace surrogate.
For $A\coloneqq \operatorname{diag}(1,\varepsilon,\varepsilon)$, $k=2$, and
$0<\varepsilon<1$, the largest determinant-root score is
$\sqrt{\varepsilon}$, whereas
$\max_{X\in\T_3}\tr(AX)=1$. Thus a linear trace score can lose an
arbitrarily large factor. Even the spectral cap used in Yuan's
dissertation \cite[Section~5.6]{YuanThesis2022} does not repair this gap:
$X\coloneqq \diag(1/2,1/2,0)$ satisfies $X\preceq I/2$ and $\tr X=1$,
but its trace score is $(1+\varepsilon)/2$.
Identity~\eqref{relax:direct-exact} is the essential property of our
relaxed score.

\subsection{The finite-rank rounding factor}\label{sec:rounding}
For a relaxed matrix of rank $r\ge k$ and $-k\beta<p\le1$, define
\begin{equation}\label{intro:bfinite}
 \directfactor{r}{k}{p}
 \coloneqq \frac rk\left[
 \prod_{j=1}^k
 \frac{\Gamma(j\beta+p/k)\Gamma((r-j+1)\beta)}
      {\Gamma(j\beta)\Gamma((r-j+1)\beta+p/k)}
 \right]^{1/p}\quad(p\ne0),
\end{equation}
with continuous value
\begin{equation}\label{round:bfinitezero}
 \directfactor{r}{k}{0}\coloneqq \frac rk
 \exp\left(\frac1k\sum_{j=1}^k
       [\psi(j\beta)-\psi((r-j+1)\beta)]\right).
\end{equation}
Its limit as $r\to\infty$ is $\factor{k}{p}$.

\begin{theorem}[Gaussian rounding]\label{round:main}
Let $A_1,\ldots,A_m\succeq0$, let $w_i>0$ sum to one, and let
$X\in\T_n$ have rank $r\ge k$. Let $Z$ be an $n\times k$
standard Gaussian matrix. For $-k\beta<p\le1$, there is a realization
of $Z$ for which the frame
$P\coloneqq G(G^*G)^{-1/2}$, with $G\coloneqq X^{1/2}Z$, satisfies
\begin{equation}\label{round:guarantee}
 \Mean_{p,w}[v_{A_i}(P)]\ge
 \directfactor{r}{k}{p}\Mean_{p,w}[f_{A_i}(X)].
\end{equation}
The factor obeys
$\factor{k}{p}\le\directfactor{r}{k}{p}\le1$.
If $r=k$, the range of $X$ gives factor one for every finite power.
\end{theorem}
The proof is in Section~\ref{app:rounding-main}. Here $X$ serves as
the covariance of each sampled column; success probabilities are
bounded later in Appendix~\ref{app:rounding}.

Cauchy--Binet writes each Gaussian numerator $\det(G^*A_iG)$ as a
nonnegative sum of squared minors. The moment of each term is the
same square-Gaussian determinant moment, after its covariance is
factored out. Orthonormalizing the sampled columns divides every
numerator by the common Gram determinant $\det(G^*G)$. At nonzero
power, a determinant-weighted probability measure cancels this
denominator; at zero, expected logarithms do so. The resulting
square and rectangular determinant moments give
\eqref{intro:bfinite}.

The change of measure is a proof device. The algorithm draws
ordinary Gaussian matrices. Section~\ref{app:rounding} supplies
a success probability that is polynomially bounded uniformly
in the covariance and the input power, away from the negative
moment boundary.

\subsection{The finite-dimensional algorithm}\label{sec:algorithms}
\begin{theorem}[Finite-dimensional approximation]\label{round:algorithm}
Fix the field and $0<\delta<1$. Let $1\le k<n$ and let the
rational input power satisfy
$-(1-\delta)k\beta\le p\le1$.
For rational PSD observations, positive rational weights, and
rational $\eps>0$, a randomized algorithm returns, with probability
at least $2/3$, a rationally encoded subspace of value at least
\[
 (\directfactor{n}{k}{p}-\eps)\OPT_{p,k}
\]
in time polynomial in the complete input bit length and
$\max(1,1/\eps)$.
\end{theorem}
The proof, including numerical precision, is in
Section~\ref{app:finite-algorithm}. The procedure is:
\begin{enumerate}
\item Approximately maximize $F$ over $\T_n$, with a small
positive covariance floor.
\item Draw independent Gaussian columns with that covariance
and keep their span.
\item Repeat polynomially many times and retain the subspace
with the largest computed objective value.
\end{enumerate}
The polynomial $\det(I+tXA_i)$ gives the required elementary
symmetric coefficient and its derivatives, without enumerating
the minors used in the analysis. The implementation does not evaluate
the Gamma formulas; they appear only in the analysis.

The finite-rank factor is strictly larger than $\factor{k}{p}$
except at $(k,p)=(1,1)$. For $k\ge2$, power-mean
monotonicity gives the explicit surplus
\[
 \frac{\directfactor{n}{k}{p}}{\factor{k}{p}}
 \ge \exp\left(\frac{k-1}{2n}\right).
\]
Together with a separate scalar estimate, this surplus pays for
optimization and output precision, yielding the exact factor in
Theorem~\ref{alg:attainment}. As $p\uparrow1$ in the scalar case,
both the finite-rank surplus $\directfactor{n}{1}{p}-\factor{1}{p}$
and the gap $1-\factor{1}{p}$ vanish
proportionally to $1-p$. Their cancellation keeps the sample
count polynomial in the input bit length.

\section{Matching computational hardness}
\label{hard:section}\label{sec:hardness}
\label{tensor:section}\label{sec:arithmetic-hardness}

The approximation factors in Section~\ref{sec:guarantees} are sharp
for each fixed selected dimension and power. The reductions below
use equal weights and have perfect completeness: in the YES case,
one subspace gives the identity compression for every observation.
All parameters in the theorem statements are fixed before the
input size grows. We write $\Mean_p$ for the equally weighted mean.

\begin{theorem}[Hardness for $p<1$]
\label{hard:main}
Fix the field, $k\ge1$, a rational $p<1$, and
$0<\varepsilon<1-\factor{k}{p}$. It is NP-hard to distinguish
rational PSD instances with the following properties:
\begin{enumerate}
\item There is $P\in\St(n,k)$ with $P^*A_iP=I_k$ for every $i$.
\item Every $P\in\St(n,k)$ satisfies
\[
 \Mean_p\bigl[v_{A_1}(P),\ldots,v_{A_m}(P)\bigr]
 \le\factor{k}{p}+\varepsilon.
\]
\end{enumerate}
The observations have rank at most $k$, and their average is an
orthogonal projection. The same gap holds when every observation
is positive definite, with the rank restriction removed and the
average a positive definite contraction.
\end{theorem}

The proof is in Appendix~\ref{app:hard:subcritical}; the rational
encoding and positive definite perturbation are in
Appendix~\ref{hard:encoding}. In particular, when
$p\le-k\beta$ the NO bound can be made arbitrarily small, excluding
every positive constant approximation factor. Geometric power
follows from the positive-power result by continuity and
power-mean monotonicity.

\begin{theorem}[Hardness at $p=1$]
\label{tensor:main}
Fix the field, $k\ge2$, and $0<\varepsilon<1-\factor{k}{1}$.
It is NP-hard to distinguish rational PSD instances with the
following properties:
\begin{enumerate}
\item There is $P\in\St(n,k)$ with $P^*A_iP=I_k$ for every $i$.
\item Every $P\in\St(n,k)$ satisfies
\[
 \frac1m\sum_{i=1}^m v_{A_i}(P)\le\factor{k}{1}+\varepsilon.
\]
\end{enumerate}
The observations have rank at most $k$, and their average is an
orthogonal projection. The same gap holds for positive definite
observations, with the rank restriction removed and the average
a positive definite contraction.
\end{theorem}

The proof is in Appendix~\ref{app:hard:arithmetic}, with the same
encoding and perturbation arguments in Appendix~\ref{hard:encoding}.
The excluded case $k=1$, $p=1$ is spectral optimization.

Both proofs begin with the block-subspace promise of
Bri\"et, Regev, and Saket \cite{BrietRegevSaket2017}, in the form
given by Ebrahimnejad, Nagda, and Oveis Gharan
\cite{EbrahimnejadNagdaOveisGharan2024}. Randomization over a finite
family of rational unitary matrices makes coefficients spread among
many coordinates obey a central limit comparison with Gaussian
determinant moments. For $p<1$,
subquadratic moments suppress the contribution of the small set of
exceptional source blocks, whose coefficients are too concentrated
for that comparison. This step fails at $p=1$:
even a small exceptional set may carry all coefficient energy.

The arithmetic reduction takes a fixed number $d$ of tensor copies
of the source. Each measurement map shares one tensor factor across
its output rows and samples the other factors independently in
each row. The proof partitions the coefficient indices using the
same good and exceptional sets for every choice of shared factor.
The resulting Cartesian products partition the coefficient energy.
A product with an exceptional independent factor contributes
little; the all-good product has the Gaussian factor $\factor{k}{1}$.
A product with one exceptional factor can retain its full energy
only when that factor is shared. Averaging the $d$ choices bounds
soundness by
$\max\{\factor{k}{1},1/d\}$, up to arbitrary fixed slack.
Choosing $1/d\le\factor{k}{1}$ gives the sharp factor.
After the independent replacements, each row is a mixture of
centered Gaussian rows. Covariance averaging controls the shared
factor with an error uniform in the source block length and all
contraction dimensions. Thus the analytic parameters can be fixed
before the source theorem supplies its block length.

The reductions are deterministic and polynomial for each fixed
$(p,k,\varepsilon)$. Their constant gadgets need not be uniform in
these parameters. Hardness at fixed selected dimension nevertheless excludes
any improved guarantee valid for arbitrary input dimensions.
The promise gaps imply deterministic approximation hardness under
$\mathrm P\ne\mathrm{NP}$ and randomized feasible-output hardness
under $\mathrm{NP}\ne\mathrm{RP}$; the stated assumption
$\mathrm{NP}\not\subseteq\mathrm{BPP}$ suffices for both feasible
outputs and bounded-error value estimation. The verification
argument, including rational representations of subspaces, is
given in Appendix~\ref{hard:encoding}.

\subsection{Hardness already holds for subspace observations}

A rank-$k$ observation separates into a scale and its range
projection. Identity~\eqref{obs:scale-identity} below also connects the hardness theorem to
the exact two-observation result in the next section.

\begin{lemma}[Projective normalization]\label{obs:normalization}
Let $A_i\succeq0$ have rank $k$, let $\Pi_i$ be its range
projector, and put $d_i\coloneqq e_k(A_i)^{1/k}>0$. For $p\ne0$, set
$S\coloneqq \sum_iw_i d_i^p$ and $\widetilde w_i\coloneqq w_i d_i^p/S$. Then
\begin{equation}\label{obs:scale-identity}
 \Mean_{p,w}[v_{A_i}(P)]
 =S^{1/p}\Mean_{p,\widetilde w}[v_{\Pi_i}(P)]
 \qquad(P\in\St(n,k)).
\end{equation}
At $p=0$, the multiplier is $\prod_i d_i^{w_i}$ and the weights
remain unchanged. For rational data and fixed rational $p$,
positive rational weights approximate the transformed weights
to any prescribed multiplicative precision in polynomial time.
\end{lemma}

The proof and the rational-weight estimate are in
Section~\ref{obs:transfer-appendix}. Perfect completeness lets the
hardness reduction reject an instance containing a deficient-rank
observation immediately. Otherwise all observations already have
rank $k$ and can be normalized directly.

\begin{theorem}[Hardness for subspace observations]\label{obs:transfer}
Fix the field, $k$, and a rational $p\le1$.
For every fixed $\rho$ with $\factor{k}{p}<\rho<1$, obtaining a
$\rho$-approximation is NP-hard when all observations are
rational rank-$k$ orthogonal projections with positive rational
weights. At $p=0$, equal weights suffice.
Here $\factor{k}{p}=0$ when $p\le-k\beta$, so every positive
constant-factor approximation is NP-hard in that range.
\end{theorem}

The proof is in Section~\ref{obs:transfer-appendix}.
This transfer concerns approximation ratios; its additional promise
properties are discussed in Section~\ref{obs:transfer-appendix}.

\section{Exactness for two \texorpdfstring{rank-$k$}{rank-k} observations}
\label{sec:observations}

For a rank-$k$ observation $A=UDU^*$, where $U\in\St(n,k)$
and $D\succ0$, put $d\coloneqq \det(D)^{1/k}=e_k(A)^{1/k}$.
Its score and relaxed score satisfy
\begin{equation}\label{obs:direct-connection}
 v_A(P)=d|\det(U^*P)|^{2/k},\qquad
 f_A(X)=kd\,\det(U^*XU)^{1/k}.
\end{equation}
Indeed, the nonzero eigenvalues of $A^{1/2}XA^{1/2}$ are
those of $D^{1/2}U^*XUD^{1/2}$. Setting $Y\coloneqq kX$ gives
the relaxed score $d\det(U^*YU)^{1/k}$.
Theorem~\ref{obs:purification} below shows that the trace relaxation
is exact for two such observations.

\subsection{Two-view purification}
\label{proof:purification}

\begin{theorem}[Two-view purification]\label{obs:purification}
Let $U_1,U_2\in\St(n,k)$. For every $Y\succeq0$ with $\tr Y=k$,
there is $P_{\mathrm{pur}}\in\St(n,k)$ such that
\begin{equation}\label{obs:purification-inequalities}
 |\det(U_i^*P_{\mathrm{pur}})|^2
 \ge\det(U_i^*YU_i),\qquad i=1,2.
\end{equation}
The construction is deterministic and uses principal vectors
and square roots of nonnegative scalars.
\end{theorem}

\begin{proof}
Changing the basis within either observation space leaves these
determinants unchanged. The principal-angle decomposition
\cite{BjorckGolub1973,Halmos1969} gives $k$ mutually orthogonal
blocks $E_j$ of dimension one or two whose orthogonal sum contains
both observation spaces. In orthonormal coordinates within these blocks, the columns
of $U_1$ and $U_2$ may therefore be taken to be
\[
 u_j\coloneqq e_j,\qquad
 v_j\coloneqq \cos\phi_j\,e_j+\sin\phi_j\,f_j,\qquad
 0\le\phi_j\le\pi/2,
\]
where $\phi_j$ is the $j$th principal angle. When $\phi_j=0$,
the block is the line spanned by $e_j$; otherwise its basis is $e_j,f_j$.
Let $C_j$ be the compression of $Y$ to $E_j$ and put
$t_j\coloneqq \tr C_j$. Orthogonality of the blocks gives
$\sum_jt_j\le\tr Y=k$; they need not span the ambient space.
If some $t_j=0$, both target determinants in
\eqref{obs:purification-inequalities} vanish, so any frame suffices.

Otherwise, in a two-dimensional block take
\[
 w_j\coloneqq \sqrt{(C_j)_{11}/t_j}\,e_j+
     \sqrt{(C_j)_{22}/t_j}\,f_j.
\]
Positivity gives
$\operatorname{Re}(C_j)_{12}
 \le\sqrt{(C_j)_{11}(C_j)_{22}}$, and hence
\[
 |u_j^*w_j|^2=\frac{u_j^*C_ju_j}{t_j},
 \qquad
 |v_j^*w_j|^2\ge\frac{v_j^*C_jv_j}{t_j}.
\]
In a one-dimensional block take $w_j\coloneqq e_j$.
The columns $w_j$ form an orthonormal frame $P_{\mathrm{pur}}$.
Because $w_j\in E_j$ and the blocks are orthogonal, both
$U_1^*P_{\mathrm{pur}}$ and $U_2^*P_{\mathrm{pur}}$ are diagonal.
Hadamard's inequality, the overlap bounds above, and
AM--GM ($\prod_jt_j\le1$) now give
\[
 \det(U_1^*YU_1)
 \le\prod_j u_j^*C_ju_j
 \le\left(\prod_jt_j\right)|\det(U_1^*P_{\mathrm{pur}})|^2
 \le|\det(U_1^*P_{\mathrm{pur}})|^2.
\]
The same argument applies to $U_2$.
Optionally, if the angle of $w_j$ from $e_j$ exceeds $\phi_j$,
replace $w_j$ by $v_j$; both overlaps increase. Thus all output
block angles may be chosen in $[0,\phi_j]$.
\qed
\end{proof}

The argument covers common directions, orthogonal directions
and vanishing compressed determinants, over either field.
Exactness uses the rank-$k$ structure of both observations.
Our hardness reductions use many observations and leave unresolved
whether hardness already holds for three.

\subsection{Exact objectives and efficient relative accuracy}

\begin{corollary}[Exactness for increasing aggregates]\label{obs:two-objectives}
Let $A_1,A_2\succeq0$ have rank $k$, with range frames $U_1,U_2$
and $d_i\coloneqq e_k(A_i)^{1/k}$. For every coordinatewise increasing
$G:[0,\infty)^2\to\mathbb R$,
\begin{align}
 &\sup_{P\in\St(n,k)}G(v_{A_1}(P),v_{A_2}(P))\notag\\
 &\quad=
 \sup_{\substack{Y\succeq0\\\tr Y=k}}
 G\left(d_1\det(U_1^*YU_1)^{1/k},
        d_2\det(U_2^*YU_2)^{1/k}\right).
 \label{obs:two-exact}
\end{align}
If $G$ is upper semicontinuous, both suprema are attained.
For $G=\Mean_{p,w}$ and every finite rational $p\le1$,
rational input admits a deterministic algorithm returning
a rationally encoded subspace of value at least
$(1-\varepsilon)\OPT_{p,k}$, for every rational
$0<\varepsilon<1$. Its time is polynomial in the complete
input bit length, including that of $p$, and in
$\log(1/\varepsilon)$.
\end{corollary}

The proof is in Section~\ref{proof:two-objectives}, with the complete
numerical argument in Section~\ref{obs:purification-precision}.

Identity~\eqref{obs:two-exact} uses only monotonicity of the aggregate;
concavity is needed for the algorithm. The algorithm has no
negative-moment restriction.

\section{Martin means of subspaces}
\label{sec:martin}

For observed frames $U_i\in\St(n,k)$, the Martin objective is
\begin{equation}\label{martin:loss}
 \mathcal L(P)\coloneqq -\sum_iw_i\log|\det(U_i^*P)|^2,
 \qquad P\in\St(n,k).
\end{equation}
The squared determinant is the product of the squared cosines
of the principal angles between the two subspaces.
A zero positive-weight overlap gives infinite loss.
This is the finite-dimensional squared-Martin objective used in
\cite{YangHospedales2016}, with origins in stochastic systems
\cite{Martin2000,DeCockDeMoor2002}.
We optimize over the whole Grassmannian, without ARMA
realization constraints.

Writing $\Pi_i\coloneqq U_iU_i^*$ gives
$\mathcal L(P)=-k\log\Mean_{0,w}[v_{\Pi_i}(P)]$.
Consequently the general approximation theorem supplies an
additive guarantee, while two-view purification yields
an explicit optimum for two observations.

\subsection{Additive approximation}

Projection observations permit a useful dimension reduction.
If $V_i$ and $Q$ are orthogonal-complement frames for $U_i$
and $P$, respectively, then
\begin{equation}\label{martin:complement-identity}
 |\det(U_i^*P)|^2=|\det(V_i^*Q)|^2.
\end{equation}
For a unitary matrix, a minor and its complementary minor have
equal modulus, by Jacobi's identity \cite{HornJohnson2013}.
Applying this to $[U_i\ V_i]^*[P\ Q]$ gives
\eqref{martin:complement-identity}, including singular overlaps.
Thus the loss is unchanged on complementing all
subspaces, and the Gaussian-rounding algorithm can work in dimension
$s\coloneqq \min(k,n-k)$.

\begin{theorem}[Dimension-free Martin approximation]\label{martin:main}
Let $1\le k<n$ and $s\coloneqq \min(k,n-k)$.
\begin{enumerate}
\item For rationally encoded observations and positive rational weights, a randomized
polynomial-time algorithm returns, with probability at least $2/3$,
\begin{equation}\label{martin:dimension-free}
 \mathcal L(P)\le\min_R\mathcal L(R)+\Delta_s,\qquad
 \Delta_s\coloneqq -s\log\factor{s}{0}
 =s\log(s\beta)-\sum_{j=1}^s\psi(j\beta).
\end{equation}
\item For fixed $k$ and $0<\varepsilon<\Delta_k$, obtaining additive
error $\Delta_k-\varepsilon$ uniformly over ambient dimension
is NP-hard, even with equal weights.
\item For rationally encoded observations, positive rational weights,
and rational $\varepsilon>0$, a randomized algorithm attains
\begin{equation}\label{martin:additive}
 \mathcal L(P)\le\min_R\mathcal L(R)+\Delta_{s,n}+\varepsilon,
\end{equation}
where
\begin{equation}\label{martin:finite-constant}
 \Delta_{s,n}\coloneqq -s\log\directfactor{n}{s}{0}
 =s\log(s/n)+\sum_{j=1}^s
 [\psi((n-j+1)\beta)-\psi(j\beta)].
\end{equation}
Its time is polynomial in the input bit length and
$\max(1,1/\varepsilon)$, and its success probability is at least $2/3$.
\end{enumerate}
\end{theorem}

The proof is in Section~\ref{proof:martin}.
The finite-dimensional statement is a guarantee of the
Gaussian-rounding algorithm in the smaller projection dimension.
The hardness assertion concerns improving $\Delta_k$ with $k$
fixed and ambient dimension growing.

\subsection{An explicit two-view mean}

Use the principal vectors $u_j,v_j$ and blocks $E_j$ from the proof
of Theorem~\ref{obs:purification}.

\begin{proposition}[Two-view Martin mean]\label{martin:two}
For two observations with principal angles
$\phi_j\in[0,\pi/2]$ and positive weights $w,1-w$,
a global minimum is attained at the subspace spanned by $e_j$ in
each common line and by
$\cos\theta_j\,e_j+\sin\theta_j\,f_j$ in each two-dimensional
principal block.
For $0<\phi_j<\pi/2$,
\begin{equation}\label{martin:angles}
 \tan\theta_j=
 \frac{2(1-w)\tan\phi_j}
 {1+\sqrt{1+4w(1-w)\tan^2\phi_j}}.
\end{equation}
At $\phi_j=0$ take $\theta_j=0$; at $\phi_j=\pi/2$ take
$\cos^2\theta_j=w$. Equal weights give principal-angle midpoints.
\end{proposition}

The proof in Section~\ref{proof:martin-two} uses two-view
purification to restrict to principal blocks and then maximizes
one strictly concave logarithmic function in each block.
Thus the same relaxation supports sharp general approximation
and an explicit global solution in the two-view case.

\section*{Statements and declarations}
\paragraph{Funding.} The author received no funding for this work.
\paragraph{Competing interests.} The author declares no competing interests.
\paragraph{Data availability.} No empirical datasets were used or generated.
All approximation factors are specified by exact formulas in the article.
\paragraph{AI acknowledgement.}
AI language models, including GPT-6 and Claude Opus, were used as supporting tools throughout the research, writing and editing process. The author retains full responsibility for the final content and conclusions.

\appendix
\section{Gaussian determinant rounding}\label{app:rounding}
We prove Theorem~\ref{round:main} and the ordinary-sampling estimates
used by the algorithms. All Gaussian normalizations are those of
Section~\ref{sec:formulation}.

\subsection{Numerator and Gram determinant moments}
For a square standard Gaussian matrix $H\in\F^{k\times k}$, put
\begin{equation}\label{round:D}
 D_{k,\beta}(u)\coloneqq \E|\det H|^{2u}
 =\beta^{-ku}\prod_{j=1}^k
       \frac{\Gamma(j\beta+u)}{\Gamma(j\beta)},\qquad u>-\beta.
\end{equation}
For $Z\in\F^{r\times k}$, $r\ge k$, the rectangular moment is
defined for $u>-(r-k+1)\beta$ by
\begin{equation}\label{round:R}
 R_{r,k,\beta}(u)\coloneqq \E\det(Z^*Z/r)^u
 =(r\beta)^{-ku}\prod_{j=1}^k
       \frac{\Gamma((r-j+1)\beta+u)}
            {\Gamma((r-j+1)\beta)}.
\end{equation}
These are Bartlett's real and complex Gram determinant formulas
\cite{EdelmanRao2005,Goodman1963Determinant}. Each Gram determinant
is a product of independent Gamma variables, with shape parameters
$(r-j+1)\beta$ and rate $\beta$.

\begin{lemma}\label{round:numerator}
Let $A,X\succeq0$, $G\coloneqq X^{1/2}Z$, and suppose
\[
 C_A(X)\coloneqq e_k(A^{1/2}XA^{1/2})>0.
\]
Then
\begin{equation}\label{round:num-moment}
 \E\det(G^*AG)^u
 \begin{cases}
 \ge D_{k,\beta}(u)C_A(X)^u,&0<u\le1,\\
 \le D_{k,\beta}(u)C_A(X)^u,&-\beta<u<0.
 \end{cases}
\end{equation}
Moreover,
\begin{equation}\label{round:num-log}
 \E\log\det(G^*AG)
 \ge\log C_A(X)+\sum_{j=1}^k\psi(j\beta)-k\log\beta.
\end{equation}
\end{lemma}
\begin{proof}
For each $k$-element coordinate set $I$, let $B_I$ be the
corresponding rows of $A^{1/2}$. Cauchy--Binet gives
\[
 \det(G^*AG)=\sum_{|I|=k}|\det(B_IG)|^2,\qquad
 C_A(X)=\sum_{|I|=k}\det(B_IXB_I^*).
\]
If $c_I\coloneqq \det(B_IXB_I^*)>0$, the term
$|\det(B_IG)|^2/c_I$ has the square Gaussian determinant law.
Terms with $c_I=0$ vanish almost surely. Thus
\[
 \det(G^*AG)=C_A(X)\sum_I\theta_ID_I,\qquad
 \theta_I\coloneqq c_I/C_A(X),\quad \sum_I\theta_I=1,
\]
where every retained $D_I$ has the law $|\det H|^2$.
Scalar Jensen proves all three assertions. Independence between
the minors is unnecessary. One positive term controls the negative
logarithmic tail, and a Gaussian norm bound controls the positive
tail, so the displayed expected logarithm is finite. \qed
\end{proof}

\begin{lemma}\label{round:denominator}
If $X\succ0$ is $r\times r$, $\tr X=1$, and $u\coloneqq p/k$, then
\begin{equation}\label{round:den-moment}
 \E\det(Z^*XZ)^u
 \begin{cases}
 \le R_{r,k,\beta}(u),&0<p\le1,\\
 \ge R_{r,k,\beta}(u),&-k\beta<p<0.
 \end{cases}
\end{equation}
Also,
\[
 \E\log\det(Z^*XZ)
 \le-k\log(r\beta)+\sum_{j=1}^k\psi((r-j+1)\beta).
\]
\end{lemma}
\begin{proof}
For fixed $Z$, the determinant power is concave in $X$ for $0<p\le1$.
For $p<0$, it is convex as the exponential of
$(p/k)\log\det(Z^*XZ)$, a convex function.
After expectation, Gaussian invariance makes the function unitarily
invariant. Its value is therefore unchanged by conjugating $X$ by
any permutation of its eigenbasis. The average of these conjugates
is $I_r/r$; Jensen bounds the value at $X$ by the value there, with
the appropriate direction for each power. Log-determinant concavity
gives the last assertion by the same averaging argument.
The moments are finite because $X\succ0$ and $u>-\beta$. \qed
\end{proof}

\subsection{Proof of the rounding theorem}\label{app:rounding-main}
\begin{proof}[of Theorem~\ref{round:main}]
Write $C_i\coloneqq C_{A_i}(X)$. If $p\le0$ and some $C_i=0$,
the assertion is immediate. For $p>0$, discard all such zero scores,
normalize the remaining weights, and restore their common factor
after rounding; if none remain, the assertion is again immediate.
We may therefore assume $C_i>0$ for every retained observation.
Work on the range of $X$ and draw
$G\coloneqq X^{1/2}Z$, $P\coloneqq G(G^*G)^{-1/2}$.
Write $N_i(G)\coloneqq \det(G^*A_iG)$. Then
\begin{equation}\label{round:ratio}
 v_{A_i}(P)^k=\frac{N_i(G)}{\det(G^*G)}.
\end{equation}
For $p\ne0$, put $u\coloneqq p/k$ and introduce the probability density
\begin{equation}\label{round:tilt}
 \omega(Z)\coloneqq \frac{\det(Z^*XZ)^u}{\E\det(Z^*XZ)^u}
\end{equation}
relative to the ordinary Gaussian law. Under the resulting law $\nu$,
\[
 \E_\nu\sum_iw_i v_{A_i}(P)^p
 =\frac{\sum_iw_i\,\E N_i(G)^u}
        {\E\det(Z^*XZ)^u}.
\]
The two lemmas bound this expression below for $p>0$, and above for
$p<0$, by
\[
 \frac{D_{k,\beta}(u)}{R_{r,k,\beta}(u)}
 \sum_iw_i C_i^u.
\]
There is therefore a realization of $Z$ for which the same bound holds.
Taking the $1/p$ power proves \eqref{round:guarantee}, with the
required reversal when $p<0$, since $f_{A_i}(X)=kC_i^{1/k}$.
At zero, take expected logarithms in \eqref{round:ratio}.

For $p>0$, Jensen applied to $Z^*Z/r$ gives
$R_{r,k,\beta}(p/k)\le1$; the inequality is reversed for $p<0$.
The logarithmic assertion follows in the same way.
Consequently $\directfactor{r}{k}{p}\ge\factor{k}{p}$.
To obtain the upper bound, apply the existence assertion just proved
to $A=I_r$ and $X=I_r/r$. Every frame has score one, whereas
$f_{I_r}(I_r/r)=(k/r)\binom rk^{1/k}$. Therefore
\begin{equation}\label{round:factor-bounds}
 \factor{k}{p}\le\directfactor{r}{k}{p}
 \le\frac{r}{k\binom rk^{1/k}}\le1.
\end{equation}
The last inequality follows by multiplying
$(r-j)/(k-j)\ge r/k$ for $0\le j<k$.
Finally, if $X=P\Lambda P^*$ has rank $k$, then
\[
 f_{A_i}(X)=k\det(\Lambda)^{1/k}v_{A_i}(P)\le v_{A_i}(P),
\]
because $\tr\Lambda=1$. This proves the final assertion for every
finite power, including zero scores. \qed
\end{proof}

\subsection{Ordinary sampling and uniform density bounds}
We now convert guarantees under the tilted law into guarantees for
ordinary Gaussian samples.
If $\nu(E)\ge\rho$ and $\E\omega^a\le C$, H\"older gives
\begin{equation}\label{round:transfer}
 \Pr(E)\ge\rho^{a/(a-1)}C^{-1/(a-1)}.
\end{equation}
We now bound these density moments uniformly in the covariance.

Let $D$ have the square Gaussian determinant law and $\mu\coloneqq \E D=k!$.
Cauchy--Binet writes $\det(Z^*XZ)=e_k(X)H_X$, where $H_X$ is
a convex combination of variables distributed as $D$. In particular
$\E H_X=\mu$. Bartlett's factorization also gives
\[
 \frac D{k!}\ \stackrel{\mathrm{law}}=\ \prod_{j=1}^kQ_j,
 \qquad Q_j\coloneqq\frac{\gamma_j}{j\beta},\quad
 \gamma_j\sim\operatorname{Gamma}(j\beta,\text{rate }1)
\]
with independent factors. Thus $Q_j$ has mean one and rate $j\beta$.
We use $\gamma_j$ for the unit-rate Bartlett factors throughout the
appendices. Put
\begin{equation}\label{round:Rsmall}
 T(t,k)\coloneqq \frac{\E D^{t/k}}{(k!)^{t/k}}
       =\prod_{j=1}^k\E Q_j^{t/k}.
\end{equation}

\begin{lemma}[Covariance-uniform density bounds]\label{round:density}
Fix $0<\delta<1$. There are constants $c=c_{\delta,\beta}>0$
and $C=C_{\delta,\beta}\ge1$, depending only on the field and
$\delta$, such that, for every covariance rank $r\ge k\ge1$
and $-(1-\delta)k\beta\le p\le1$,
\[
 c\le\factor{k}{p}\le\directfactor{r}{k}{p}\le1.
\]
The lower constant $c$ may be chosen rational. For the density
\eqref{round:tilt},
\begin{align}
 \E\omega^2&\le3c^{-2},&&p>0,\label{round:positive-density}\\
 \E\omega^{1+\delta/2}&\le Ck,&&p<0.\label{round:negative-density}
\end{align}
\end{lemma}
\begin{proof}
Each $Q_j$ is distributed as the average of $j$ independent
copies of $Q_1$. For $u<0$, Jensen gives
$\E Q_j^u\le\E Q_1^u$. Hence, for negative allowed $p$,
\[
 \factor{k}{p}
 =\frac{(k!)^{1/k}}k\,T(p,k)^{1/p}
 \ge e^{-1}\Mean_{p/k}[Q_1]
 \ge e^{-1}\Mean_{-(1-\delta)\beta}[Q_1]>0.
\]
Here $(k!)^{1/k}/k\ge e^{-1}$ follows from
$\log(k!)\ge\int_1^k\log x\dd x=k\log k-k+1$.
Power-mean monotonicity and continuity extend
the same bound to $p\ge0$. Choose a positive rational $c$ below
this fixed constant. The upper bounds follow from
\eqref{round:factor-bounds}.

For negative $p$, set $a\coloneqq 1+\delta/2$, $v\coloneqq -ap/k$, and
$v_*\coloneqq (1-\delta/2)\beta$. Then $0<v\le v_*<\beta\le1$.
Jensen for the mixture, followed by $\E H_X^{p/k}\ge\mu^{p/k}$,
gives
\[
 \E\omega^a\le T(ap,k)=\prod_{j=1}^k\E Q_j^{-v}.
\]
For $j\beta>1$, concavity of the $v$th power and the Gamma
recurrence imply
\[
 \E Q_j^{-v}\le(\E Q_j^{-1})^v
              =\left(\frac{j\beta}{j\beta-1}\right)^v.
\]
Write $j_0\coloneqq 1/\beta\in\{1,2\}$. For $k\ge j_0$ the product of the
ratios for $j=j_0+1,\ldots,k$ is $\binom{k}{j_0}\le k^{j_0}$.
The remaining factors satisfy
$\E Q_j^{-v}\le\E Q_j^{-v_*}$ by H\"older interpolation and
Jensen. Thus, with
\[
 C\coloneqq \prod_{j=1}^{j_0}\E Q_j^{-v_*},
\]
the product is at most $Ck^{j_0v}\le Ck$.
For $k<j_0$, the same bound holds by omitting the absent factors,
each of which is at least one. This proves
\eqref{round:negative-density}.

For positive $p$, write $u\coloneqq p/k$. If $2u\le1$, Jensen and the
numerator mixture estimate give
\[
 \E\omega^2
 \le\frac{\mu^{2u}}{(\E D^u)^2}=T(p,k)^{-2}\le c^{-2}.
\]
Indeed $\factor{k}{p}=[(k!)^{1/k}/k]T(p,k)^{1/p}$,
so $T(p,k)\ge c^p\ge c$.
The remaining case is $k=1$, $1/2<p\le1$.
Convexity for the mixture and interpolation between its first
two moments give
$\E H_X^{2p}\le\E D^{2p}\le(\E D^2)^{2p-1}\le3$,
whereas $\E H_X^p\ge\E D^p\ge c$. This proves
\eqref{round:positive-density}. \qed
\end{proof}

\subsection{Success probability, including powers near zero}
Assume all relaxed scores $f_{A_i}(X)$ are positive and $X$ maximizes the
relaxed objective, with value $S^\star>0$. Set
\[
 Y\coloneqq \frac{\Mean_{p,w}[v_{A_i}(P)]}{S^\star},\qquad
 b\coloneqq \directfactor{r}{k}{p},\qquad L\coloneqq -\log b.
\]
Exactness at subspaces gives $0\le Y\le1$.
The rounding proof yields
$\E_\nu Y^p\ge b^p$ for $p>0$,
$\E_\nu Y^p\le b^p$ for $p<0$, and
$\E\log Y\ge\log b$ at zero.
For $s>0$, apply these inequalities to the event
$E\coloneqq \{Y\ge be^{-s}\}$. Bounding $Y$ by one on $E$ gives
\begin{equation}\label{round:success}
 \nu(E)\ge
 \begin{cases}
 \displaystyle
 \frac{e^{ps}-1}{e^{p(L+s)}-1}
 \ge e^{-(L+s)}\frac{s}{L+s},&p>0,\\[6pt]
 \displaystyle
 \frac{1-e^{-|p|s}}{1-e^{-|p|(L+s)}}
 \ge\frac{s}{L+s},&p<0,\\[6pt]
 \displaystyle\frac{s}{L+s},&p=0.
 \end{cases}
\end{equation}
At zero, $\nu$ denotes the ordinary law. The first lower bound
uses $0<p\le1$; the second follows from monotonicity of
$(1-e^{-x})/x$. Retaining the denominators makes the estimate
uniform as $p$ tends to zero.

By Lemma~\ref{round:density}, $0\le L\le-\log c$.
With $s\coloneqq -\log(1-\varepsilon)$ and $0<\varepsilon\le1/2$,
the tilted success probability is therefore at least
$c'\varepsilon$, for a fixed $c'>0$.
Combining \eqref{round:transfer} with the density lemma gives,
for negative powers,
\begin{equation}\label{round:ordinary-prob}
 \Pr\{Y\ge(1-\varepsilon)b\}
 \ge(c'\varepsilon)^{1+2/\delta}(Ck)^{-2/\delta}.
\end{equation}
For positive powers the lower bound is
$(c'\varepsilon)^2c^2/3$, and at zero it is $c'\varepsilon$.
Thus ordinary independent Gaussian samples suffice with a
polynomial count, uniformly in $k$, the input power, and the
covariance. Section~\ref{app:implementation} extends the argument
to approximate relaxed optima and rational samples.

\section{Algorithms and numerical precision}
\label{app:algorithms}\label{app:implementation}
We prove Theorems~\ref{round:algorithm} and~\ref{alg:attainment}.
All polynomial bounds may depend on the field and the fixed margin
$\delta$. We use weak convex optimization with evaluation and
separation oracles \cite{GroetschelLovaszSchrijver1988}.

\subsection{Zero scores and a positive optimization scale}
If $\rank A_i<k$, the corresponding volume is identically zero.
For $p\le0$, any such positive-weight score makes the optimum zero.
For $p>0$, delete these zero scores and normalize the remaining weights.
If their original sum is $W>0$, the original objective is $W^{1/p}$
times the normalized objective; if $W=0$, the problem is trivial.
The factor $W^{1/p}$ is retained symbolically and cancels from every
approximation ratio. Its logarithm need not have polynomial magnitude
when $p$ is extremely small.

Assume henceforth that every retained $A_i$ has rank at least $k$
and $k<n$. Write $F(X)\coloneqq \Mean_{p,w}[f_{A_i}(X)]$ and restrict to
\[
 \T_n^\theta
 \coloneqq \{X\succeq(\theta/n)I:\tr X=1\},\qquad 0<\theta\le1/4.
\]
For any $X\in\T_n$, the mixture
$X_\theta\coloneqq (1-\theta)X+\theta I/n$ is in $\T_n^\theta$ and
$F(X_\theta)\ge(1-\theta)F(X)$.
Moreover, for $X\in\T_n^\theta$,
\begin{equation}\label{app:positive-scale}
 e_k(A_i^{1/2}XA_i^{1/2})\ge(\theta/n)^k e_k(A_i)>0.
\end{equation}
A positive rational principal minor gives a polynomial-bit lower
bound for $e_k(A_i)$; matrix norm bounds give an upper bound.
Thus every relaxed score, and their mean, has a logarithm bounded by
a polynomial in the input length and $\log(1/\theta)$.
In particular $S_0\coloneqq F(I/n)>0$ provides an accessible relative scale.

\begin{lemma}[Logarithmic stability and derivative bounds]\label{alg:stability}
For every finite real $p$, positive weights, and positive vectors $a,b$,
\begin{equation}\label{app:log-stability}
 |\log \Mean_{p,w}[a]-\log \Mean_{p,w}[b]|
 \le\|\log a-\log b\|_\infty.
\end{equation}
If a differentiable function $F$ on the positive definite cone is
Loewner-increasing and homogeneous of degree $\rho>0$, then, on
$X\succeq aI$,
\begin{equation}\label{app:gradient}
 \nabla F(X)\succeq0,\qquad
 \tr\nabla F(X)\le\frac{\rho F(X)}a.
\end{equation}
\end{lemma}
\begin{proof}
The derivative of $\log\Mean_{p,w}$ with respect to $\log a_i$ is
$w_i a_i^p/\sum_jw_ja_j^p$ for $p\ne0$, and $w_i$ at zero.
These derivatives form a probability vector. Integration along a segment
in logarithmic coordinates proves \eqref{app:log-stability}.
Monotonicity gives the PSD gradient, and Euler's identity
$\tr(X\nabla F)=\rho F$ gives \eqref{app:gradient}. \qed
\end{proof}

Our objective has degree one. The lemma bounds its first derivatives
using the spectral floor and objective bounds, independently of the
magnitude of the power. The needed coefficient
$e_k(A_i^{1/2}XA_i^{1/2})=[t^k]\det(I+tXA_i)$ is computed by
polynomial arithmetic; here $[t^k]$ denotes the coefficient of $t^k$.
Interpolation at $n+1$ rational nodes recovers every coefficient
with polynomial bit complexity. Directional derivatives are obtained from
\[
 \left.\frac{d}{dh}\det(I+t(X+hH)A_i)\right|_{h=0}
 =t\,\tr\bigl(\operatorname{adj}(I+tXA_i)HA_i\bigr).
\]
Adjugates and interpolation compute these coefficients without
assuming invertibility. The positive lower bounds determine the
precision needed for roots, logarithms, and their derivatives.

The constraints defining $\T_n^\theta$ admit spectral separation.
In the affine space $\tr X=1$, the set contains a ball about $I/n$
of radius $(1-\theta)/(2n)$ and has bounded diameter.
Weak optimization with the preceding oracles therefore produces
rational $X\in\T_n^\theta$ with
\[
 F(X)\ge\max_{\T_n^\theta}F-\theta S_0
       \ge(1-2\theta)\max_{\T_n}F.
\]
Certified interval evaluation of $S_0$ supplies a positive rational
lower scale. The deterministic optimization cost is polynomial in
the input length, $\log(1/\theta)$, and the logarithm of the inverse
requested relative accuracy. Using smaller internal tolerances
permits a final mixture with $I/n$ before rational rounding,
preserving the trace and the positive floor within the stated
error budget.
The logarithmic dependence on accuracy is needed below when
$1-p$ is exponentially small in its encoding length.

\subsection{Evaluating rational powers}
For positive relaxed scores put $\ell_i\coloneqq \log f_i$ and evaluate
\begin{equation}\label{app:log-mean}
 \log\Mean_{p,w}[f]
 =a+\frac1p\log\left(\sum_iw_i e^{p(\ell_i-a)}\right).
\end{equation}
At zero use $\sum_iw_i\ell_i$.
Near zero, a centered series, or extra precision for the logarithm
near one, avoids cancellation. The additional precision is
proportional to $\log(1/|p|)$, bounded by the input bit length.
For negative $p$, choose $a\coloneqq \min_i\ell_i$; the exponential arguments
are nonpositive and their weighted sum lies between $w_{\min}$ and
one. Terms below the requested tolerance may be discarded.
Even if $|p|$ is large, only $O(\log|p|)$ additional operand bits
are needed. The same centering evaluates the probability weights
in the gradient. Since $p\le1$, positive exponents introduce no
additional dependence on a large power.
These procedures, together with the preceding polynomial
logarithmic bounds, give evaluation and gradient oracles to
any specified polynomial-bit accuracy.

\subsection{Success from an approximate optimum}
Let $S^\star\coloneqq \max_{\T_n}F$ and
$q\coloneqq F(X)/S^\star\ge1-2\theta$.
For the normalized rounded value
\[
 Y\coloneqq \Mean_{p,w}[v_{A_i}(P)]/S^\star
\]
one still has $0\le Y\le1$. The moment estimates in
\eqref{round:success} use $bq$ in place of
$b\coloneqq \directfactor{n}{k}{p}$. Since
$-\log(bq)\le-\log c_{\delta,\beta}-\log(1-2\theta)$,
the same inverse-polynomial success bounds apply.
Choose optimization and sampling losses to be sufficiently small
constant multiples of the requested factor slack.
Independent repetition gives a candidate of value at least
$(\directfactor{n}{k}{p}-\eps/2)S^\star$ with probability at
least $5/6$. The algorithm does not use the unknown $S^\star$;
it compares objective values and retains the largest.

\subsection{Small-determinant bounds and numerical precision}
For square standard Gaussian $H\in\F^{k\times k}$, Bartlett gives
$|\det H|^2=\beta^{-k}\prod_{j=1}^k\gamma_j$ in distribution,
with independent unit-rate factors
$\gamma_j\sim\operatorname{Gamma}(j\beta,\text{rate }1)$.
Each $\gamma_j$ stochastically dominates $\gamma_1$. Consequently
\begin{equation}\label{app:small-det}
 \Pr\{|\det H|^2<t^k\}
 \le C_\beta k\,t^\beta,\qquad
 C_\beta\coloneqq \frac{\beta^\beta}{\Gamma(\beta+1)},\quad 0<t\le1.
\end{equation}
Indeed, a product below $t^k$ has some $\gamma_j<\beta t$;
integrate the Gamma density bound for that factor and take a union bound.

For observation $A_i$, the Cauchy--Binet expansion in
Lemma~\ref{round:numerator} has at most $J\coloneqq \binom nk\le2^n$
terms and total covariance coefficient
$C_i\coloneqq e_k(A_i^{1/2}XA_i^{1/2})$.
One coefficient is at least $C_i/J$.
Marginally, its squared minor is that coefficient times
$|\det H|^2$. Therefore
\[
 \Pr\{\det(G^*A_iG)<(C_i/J)t^k\}\le C_\beta k\,t^\beta.
\]
No term needs to be found by the algorithm.
The floor and \eqref{app:positive-scale} give a lower bound for
$C_i$ with polynomial-bit logarithm.
The same reasoning applies to the Gram denominator by taking
$A_i=I_n$.

A union bound over all observations and all samples, together
with Gaussian entry tails, gives an event of probability at
least $1-\zeta$ on which every relevant determinant exceeds
$2^{-B}$ and all matrix norms are below $2^B$, where $B$ is
polynomial in the input length, $\log(1/\theta)$,
$\log(1/\zeta)$, and the logarithm of the sample count.
Independence between observation determinants is unnecessary.
For the Gram matrix, also use
\[
 \lambda_{\min}(G^*G)
 \ge\frac{\det(G^*G)}{\|G^*G\|_{\mathrm{op}}^{k-1}}.
\]
Determinant multilinearity and these conditioning bounds show
that polynomial-bit approximations to $X^{1/2}$ and to the
Gaussian entries preserve rank and approximate every positive
volume to any specified relative tolerance, with precision
polynomial in the logarithm of its inverse.
A truncated inverse Gaussian distribution function applied to
finite uniform binary strings supplies such rational samples.
Outside the tail event already included in the failure budget
above, that inverse has an exponentially bounded derivative.
Polynomially many uniform bits therefore suffice to approximate
each entry to the required absolute tolerance, under a coupling
with an exact Gaussian sample.
Certified square-root approximation on the spectral interval
$[\theta/n,1]$ has the same polynomial precision bound.

Lemma~\ref{alg:stability} transfers the relative score
accuracy to the objective uniformly in the power.
A rational output basis $R$ is evaluated by
\begin{equation}\label{app:rational-value}
 v_{A_i}(\operatorname{range}R)^k
 =\frac{\det(R^*A_iR)}{\det(R^*R)}.
\end{equation}
Both determinants are rational. Thus the algorithm needs no
exact irrational orthonormalization.
The same conditioning bounds give relative score evaluation.
Selecting the largest approximate value loses at most twice
the evaluation tolerance in logarithmic scale.

\subsection{Finite-dimensional approximation}\label{app:finite-algorithm}
\begin{proof}[of Theorem~\ref{round:algorithm}]
Cases with identically zero optimum were handled above.
We may replace $\eps$ internally by
$\min(\eps,c_{\delta,\beta}/4)$; this costs only a fixed factor
in the stated complexity and makes all intermediate guarantees positive.
Regularization makes the relaxed covariance full rank, so
Theorem~\ref{round:main} supplies $\directfactor{n}{k}{p}$.
Approximately maximize $F$, draw the polynomial number of
ordinary Gaussian samples supplied by \eqref{round:ordinary-prob}
and its nonnegative-power counterpart, and retain the best
computed subspace. Allocate half of the factor slack to
optimization and sampling, and half to discretization and
selection. Set the total coupling failure to at most $1/6$.
The ideal sampling failure is at most $1/6$, so the success
probability is at least $2/3$.
All precisions are polynomial in the complete input length
and $\log(1/\eps)$; the sample count is polynomial in the input
length and $\max(1,1/\eps)$. Since $S^\star\ge\OPT_{p,k}$,
the claimed guarantee follows. \qed
\end{proof}

\subsection{The strict finite-dimensional margin}\label{app:finite-margin}
For $k\ge2$, let $D\coloneqq \det(Z^*Z/n)$ with $Z\in\F^{n\times k}$,
and put
\[
 \mathfrak g_{n,k}
 \coloneqq \left(\prod_{j=1}^k\frac{n-j+1}{n}\right)^{1/k}.
\]
The rectangular moment formula gives $\E D=\mathfrak g_{n,k}^k$
and, including the continuous interpretation at zero,
\begin{equation}\label{alg:finite-margin}
 \frac{\directfactor{n}{k}{p}}{\factor{k}{p}}
 =\frac1{\Mean_p[D^{1/k}]}.
\end{equation}
Since $p\le1\le k$, power-mean monotonicity gives
$\Mean_p[D^{1/k}]\le(\E D)^{1/k}=\mathfrak g_{n,k}$.
AM--GM for $(n-j+1)/n$ now yields
\begin{equation}\label{alg:dimension-margin}
 \frac{\directfactor{n}{k}{p}}{\factor{k}{p}}
 \ge\frac1{\mathfrak g_{n,k}}
 \ge\frac1{1-(k-1)/(2n)}
 \ge\exp\left(\frac{k-1}{2n}\right).
\end{equation}
For $k=1$ this argument has no strict surplus; a short Gamma
estimate supplies the missing quantitative information.

\begin{lemma}[Gamma log-convexity estimate]\label{alg:gamma-jensen}
For $x>0$, $t<1$, and $x+t>0$,
\begin{equation}\label{alg:gamma-margin}
 \log x-\frac{\log\Gamma(x+t)-\log\Gamma(x)}t
 \ge\frac{1-t}{2(x+1)},
\end{equation}
with the continuous interpretation at $t=0$.
\end{lemma}
\begin{proof}
Put $g(t)\coloneqq \log\Gamma(x+t)$.
The trigamma series gives
$\psi'(y)=\sum_{j=0}^\infty(y+j)^{-2}>1/y$.
On $[\min(t,0),1]$, therefore, $g''>1/(x+1)$,
so $g(t)-t^2/[2(x+1)]$ is convex.
Interpolate between zero and one when $0<t<1$, or extrapolate
to $t<0$, and use $g(1)-g(0)=\log x$.
Division by $t$, with its appropriate sign, gives
\eqref{alg:gamma-margin}; take a limit at zero. \qed
\end{proof}
Taking $x\coloneqq n\beta$ and $t\coloneqq p$, and using $e^u-1\ge u$, gives
\begin{equation}\label{alg:scalar-margin}
 \directfactor{n}{1}{p}-\factor{1}{p}
 \ge\frac{\factor{1}{p}(1-p)}{n(2\beta+1)}
 \qquad(n\ge2,\ -\beta<p<1).
\end{equation}
Here $2(n\beta+1)\le n(2\beta+1)$.

\subsection{Scalar attainment uniformly near power one}
\label{app:scalar-attainment}
It remains to justify the last nontrivial case in
Theorem~\ref{alg:attainment}: $k=1$, $n\ge2$, and $1/2\le p<1$.
Put
\[
 s\coloneqq 1-p,\qquad b\coloneqq \factor{1}{p},\qquad
 a\coloneqq b^p,\qquad A\coloneqq \directfactor{n}{1}{p}^p.
\]
For a mean-one Gamma variable $Q$ of shape $\beta$, $a=\E Q^p$.
Power-mean monotonicity gives
$b\ge\factor{1}{1/2}=2/\pi$ over $\R$ and $\pi/4$ over $\C$.
In particular $a,b>1/2$.
Convexity of $t\mapsto\E Q^t$ and its tangent at one imply
\begin{equation}\label{app:scalar-deficit}
 1-a\le s\,\E(Q\log Q)\le s\,\E[Q(Q-1)]
       =s/\beta\le2s,
\end{equation}
where $\log Q\le Q-1$ was used.
For $x\coloneqq n\beta$, write
$r_n(p)\coloneqq x^{-p}\Gamma(x+p)/\Gamma(x)$, so $A=a/r_n(p)$.
Put $d\coloneqq s/(24n)$. Lemma~\ref{alg:gamma-jensen} gives
\begin{equation}\label{app:scalar-gain}
 A-a\ge\frac{ap(1-p)}{2(n\beta+1)}
       \ge\frac{s}{24n}=d.
\end{equation}
We used weaker bounds $a\ge1/4$, $p\ge1/2$, and $\beta\le1$
in the last estimate.
By \eqref{round:factor-bounds}, $A\le1$, so $a+d\le1$.

Obtain positive definite rational $X$ with
$q\coloneqq F(X)/S^\star\ge1-d/4$, using a covariance floor proportional
to $d/n$. The optimization precision and iteration count above
depend polynomially on
$\log(1/d)=O(\log n+\log(1/(1-p)))$, hence on the rational
input bit length. For $Y\coloneqq \Mean_{p,w}[v_{A_i}(P)]/S^\star$
under the determinant tilt, the finite-rank moment estimate gives
\[
 \E_\nu Y^p\ge Aq^p\ge A(1-d/4)\ge a+3d/4.
\]
Define $T$ only for analysis by $T\coloneqq(a+d/4)^{1/p}$.
If $E\coloneqq \{Y\ge T\}$, the bound $0\le Y\le1$ gives
\begin{equation}\label{app:scalar-success}
 \nu(E)\ge\frac{d/2}{1-a-d/4}
          \ge\frac{d}{2(1-a)}
          \ge\frac1{96n}.
\end{equation}
The factor $1-p$ has cancelled.
On $[b,1]$, the derivative of $x^p$ is at most
$b^{-1/2}\le2$, whence
\begin{equation}\label{app:scalar-buffer}
 T-b\ge d/8.
\end{equation}

For the ordinary law, diagonalize $X$ and write the Gram scalar
as $D\coloneqq \sum_j\lambda_jQ^{(j)}$, where $\lambda_j\ge0$ sum to one
and the $Q^{(j)}$ are independent copies of $Q$.
Concavity gives $\E D^p\ge a\ge1/2$.
Since $1\le2p\le2$, convexity and log-convexity of moments give
\[
 \E D^{2p}\le\E Q^{2p}
 \le(\E Q^2)^{2p-1}\le1+1/\beta\le3.
\]
Thus the tilt $W\coloneqq D^p/\E D^p$ has $\E W^2\le12$, and
\[
 \Pr(E)\ge\frac{\nu(E)^2}{12}
          \ge\frac1{12(96n)^2}.
\]
At least $12(96n)^2\log6$ ordinary samples make failure
probability at most $1/6$, independently of $1/(1-p)$.

Use the discretization argument above with relative score,
basis, and evaluation errors bounded by sufficiently small
fixed multiples of $d$, for example $d/128$ each.
The floor, small-minor bound, and sample count require only
polynomially many bits in $\log(1/d)$ and the input size.
The buffer \eqref{app:scalar-buffer} absorbs the combined errors
in the sampled candidate and selection by approximate values.
Take total coupling failure at most $1/6$.
With probability at least $2/3$, the retained rational line has
value at least $bS^\star\ge b\OPT_{p,1}$.
The constants $a,A,b,T$ are not evaluated by the algorithm.

\subsection{Attainment and assembly of the classification}\label{app:attainment}
\begin{proof}[of Theorem~\ref{alg:attainment}]
Lemma~\ref{round:density} gives a positive rational lower bound
$c_{\delta,\beta}$ for $\factor{k}{p}$, independent of the input.
For $k\ge2$, \eqref{alg:dimension-margin} gives an additive
finite-dimensional margin at least $c_{\delta,\beta}/(2n)$.
Apply Theorem~\ref{round:algorithm} with slack
$c_{\delta,\beta}/(4n)$. This is polynomial in the complete
input length and attains $\factor{k}{p}$.
When $k=n$, return the unique subspace directly.
For $k=1$ and $p\le1/2$, \eqref{alg:scalar-margin} permits
slack $c_{\delta,\beta}/[4n(2\beta+1)]$, again with polynomial
running time.
For $1/2\le p<1$, Section~\ref{app:scalar-attainment} keeps the
sample count $O(n^2)$ and requires precision only polynomial in
$\log n+\log(1/(1-p))$.
At $k=p=1$, the objective is the Rayleigh quotient of
$\sum_iw_iA_i$. An approximation to its top eigenspace gives
every prescribed factor below one in the stated running time.
Exact eigenvectors are algebraic. Rational exactness already fails
for $\left(\begin{smallmatrix}2&1\\1&1\end{smallmatrix}\right)$:
its top eigenvector has slope $(\sqrt5-1)/2$, so no nonzero
rational vector spans that line.
\qed
\end{proof}

\begin{proof}[of Theorem~\ref{intro:classification}]
For each fixed $p>-k\beta$, choose $0<\delta<1$ with
$p\ge-(1-\delta)k\beta$ and apply Theorem~\ref{alg:attainment}.
Theorems~\ref{hard:main} and~\ref{tensor:main} give hardness
for $p<1$ and $p=1$, including the positive-definite
families. Theorem~\ref{obs:transfer} gives the projection subclass.
At and below $-k\beta$ the threshold in Theorem~\ref{hard:main}
is zero. The exception at $(k,p)=(1,1)$ is the eigenvalue
problem just considered. \qed
\end{proof}

\section{Matching computational hardness}\label{app:hardness}

This appendix proves Theorems~\ref{hard:main}, \ref{tensor:main},
and~\ref{obs:transfer}, together with the normalization identity in
Lemma~\ref{obs:normalization}. We construct instances with optimum one
in the YES case and with every feasible value at most the Gaussian
threshold plus an arbitrarily small fixed slack in the NO case.
Both constructions start from a block-subspace promise. In the NO
case, it ensures that most coefficient blocks have no large individual
coefficient, allowing comparison with Gaussian sums.

A small exceptional set of blocks may nevertheless carry all the
coefficient energy. For $p<1$, moment estimates control its effect:
positive subquadratic moments bound its contribution, while negative
moments allow us to retain only the good blocks. At $p=1$, the
subquadratic estimate fails. We instead use several tensor factors
and average over which factor is shared by the output rows. This
averaging limits the contribution of concentrated energy.

Section~\ref{hard:source} supplies the common ingredients;
Sections~\ref{app:hard:subcritical} and~\ref{app:hard:arithmetic}
prove the two PSD gaps. Section~\ref{hard:encoding} handles rational
encoding, positive-definite perturbations, and the remaining powers.
Section~\ref{obs:transfer-appendix} transfers the hardness to projection
observations. Throughout, $\beta\coloneqq 1/2$ over $\mathbb R$ and $\beta\coloneqq 1$
over $\mathbb C$. All reduction parameters are fixed before the source
input size grows.

\subsection{The shared source and moment family}
\label{hard:source}

We use the following form of the block-subspace hardness theorem of
Bri\"et, Regev, and Saket \cite{BrietRegevSaket2017}, stated explicitly
as Theorem~4.3 in the full version of
Ebrahimnejad, Nagda, and Oveis Gharan
\cite{EbrahimnejadNagdaOveisGharan2024}. All norms in this statement
are ordinary, unnormalized Euclidean norms. The YES witness selects
one coordinate per block. The NO promise limits the number of blocks
that have both bounded norm and a large coordinate.

\begin{proposition}[Block-subspace promise]
\label{hard:block-source}
For any fixed $0<\tau,\delta_0,\eta<1$, there is a fixed block length $b_{\rm src}$
such that the following problem is NP-hard. Given a rational
orthogonal projection onto $\mathcal W\subseteq\mathbb R^{N b_{\rm src}}$, distinguish:
\begin{enumerate}
\item $\mathcal W$ contains a vector $v$ such that every block
$v_i\in\mathbb R^{b_{\rm src}}$ is a coordinate unit vector.
\item For every $v\in \mathcal W$ with $\|v\|^2=N$, at most $\tau N$ blocks
satisfy both $\|v_i\|\le\eta^{-1}$ and
$\|v_i\|_\infty\ge\delta_0$.
\end{enumerate}
\end{proposition}

We use only this small-parameter regime. Its cases are disjoint:
every block of the YES witness has norm and maximum coordinate one,
so all $N$ blocks satisfy the two NO tests, whereas $\tau N<N$.
Decreasing any of the three source parameters strengthens the NO
guarantee. All choices below may therefore be capped at $1/2$.

The rational representation follows directly from the original
construction. Its subspace is defined by homogeneous equalities
between coordinate sums
\cite[Section~3.1, equation~(3.1)]{BrietRegevSaket2017}.
After selecting independent rows of the resulting integer constraint
matrix $D$, the projection is
$\Pi_{\mathcal W}=I-D^T(DD^T)^{-1}D$. It is computable with polynomial bit length;
the encoding details are given in Appendix~\ref{hard:encoding}.
We use this real source over both fields. For $\F=\mathbb C$, $\mathcal W$
is complexified; its real NO promise will be applied to concatenated
real and imaginary coefficient columns.

\subsubsection{A rational unitary family with independent row phases}
\label{hard:gadget}

Let $C_t$ be a uniformly chosen cyclic permutation of $k$ coordinates,
and let $\varepsilon_1,\ldots,\varepsilon_k$ be independent uniform
signs over $\mathbb R$, or fourth-root phases over $\mathbb C$.
Write $\pi_t$ for the coordinate permutation determined by
$e_r^TC_t=e_{\pi_t(r)}^T$.
Our finite indexed family $\mathcal U_k$ consists of
\[
 U\coloneqq \diag(\varepsilon_1,\ldots,\varepsilon_k)C_t.
\]
Every outcome is a rational orthogonal or unitary matrix, and
\begin{equation}
 \E U=0,\qquad
 \E U_{ab}\overline{U_{cd}}=\frac{\delta_{ac}\delta_{bd}}{k}.
 \label{hard:isotropy}
\end{equation}
Different row phases have zero cross-moment; within one row the
uniform shift visits each column with probability $1/k$.
Over $\mathbb C$, $\E U_{ab}U_{cd}=0$ as well because each phase
has zero second moment. The support has $k2^k$ real or $k4^k$
complex parameter tuples. Its exponential dependence on $k$ is
harmless here: the selected dimension is fixed before the source
input grows. All tuples are enumerated with equal weight.
The independent row phases will also be used in the arithmetic proof.

\subsubsection{A norm-exact source-coefficient bound}

Let $\mathcal W_{\F}$ denote the real source space or its
complexification. All Euclidean and Frobenius norms are unnormalized.
The next lemma extends the scalar source promise to matrix coefficients,
uniformly in the number of columns. This uniformity will allow us to
apply it to large tensor reshapes. The proof uses orthogonal columns
so that every random-sign combination has squared norm exactly $N$,
as required by the source promise.

\begin{lemma}[Few exceptional source blocks]
\label{tensor:hilbert}
Fix $M\ge1$ and $0<d_0<1$. Suppose the source NO promise holds with
parameters $\tau_0$, $\delta_0\le d_0/2$, and
$\eta_0\le(4M)^{-1}$. If the columns of $V$ belong to
$\mathcal W_{\F}$ and $\|V\|_{\mathrm F}^2=N$, at most
$8\tau_0N$ source blocks satisfy
\[
 \|V_i\|_{\mathrm F}\le M,\qquad
 \max_a\|V_{i,a}\|_2\ge d_0.
\]
Here $V_{i,a}$ is the row at coordinate $(i,a)$.
Outside at most $(8\tau_0+M^{-2})N$ blocks, therefore, both
$\|V_i\|_{\mathrm F}\le M$ and
$\max_a\|V_{i,a}\|_2<d_0$ hold.
The estimate is independent of the number of columns.
\end{lemma}
\begin{proof}
First suppose $V$ is real. Choose an orthogonal $Q$ that
diagonalizes $V^TV$, so the columns $b_j$ of $B\coloneqq VQ$ are
orthogonal. Their row and block squared norms are unchanged.
For independent signs $\varepsilon_j$, put $x\coloneqq \sum_j\varepsilon_jb_j$.
Every outcome belongs to $\mathcal W$ and has exactly
\[
 \|x\|^2=\sum_j\|b_j\|^2=N.
\]
Fix a block in the asserted exceptional set and a row with
squared norm $\sigma^2\ge d_0^2$. The corresponding scalar
$Z\coloneqq x_{i,a}$ satisfies $\E Z^2=\sigma^2$ and
\[
 \E Z^4=3\sigma^4-2\sum_jB_{i,a,j}^4\le3\sigma^4.
\]
If $E\coloneqq \{|Z|\ge\sigma/2\}$, Cauchy--Schwarz gives
$\sigma^2\le\sigma^2/4+(\E Z^4\Pr(E))^{1/2}$, hence
$\Pr\{|Z|\ge d_0/2\}\ge3/16$.
Moreover $\E\|x_i\|^2\le M^2$, so
$\Pr\{\|x_i\|>4M\}\le1/16$.
Both source tests therefore hold with probability at least $1/8$.
The source promise bounds the expected total number of such blocks
by $\tau_0N$, proving the first assertion. At most $N/M^2$
further blocks have $\|V_i\|_{\mathrm F}>M$.
For complex $V$, apply the real argument to
$[\Re V\ \Im V]$, which preserves all row, block, and total norms.
The orthogonalization is only a proof device; it is not encoded in
the reduction. \qed
\end{proof}

\subsubsection{Gaussian comparison when every coefficient is small}
\label{hard:clt}

For a square matrix $T\in\F^{k\times k}$, put
\[
 \mathcal H(T)\coloneqq \sqrt{k}\,|\det T|^{1/k}.
\]
Let $U_j$ be independent samples from $\mathcal U_k$.
Suppose the squared Frobenius norms of the coefficient matrices $Z_j$
sum to one. Unitarity and the second-moment identity give
\begin{equation}
 \mathcal H\Bigl(\sum_jU_jZ_j\Bigr)^2
 \le\Bigl\|\sum_jU_jZ_j\Bigr\|_{\mathrm F}^2,
 \qquad
 \E\Bigl\|\sum_jU_jZ_j\Bigr\|_{\mathrm F}^2=1.
 \label{hard:quadratic-control}
\end{equation}

The next lemma improves this second-moment bound to the Gaussian
benchmark whenever every normalized coefficient is small. Its
coefficient threshold is uniform in the number of summands.

\begin{lemma}[Uniform determinant moments]
\label{hard:moment-lemma}
Fix $k\ge1$. Let $-k\beta<p<1$, $p\ne0$, $q\coloneqq 2p$, and
$\kappa>\sqrt{\factor{k}{p}}$. There is $\delta>0$ such that, for every
finite family $(Z_j)$ satisfying
\[
 \sum_j\|Z_j\|_{\mathrm F}^2=1,
 \qquad \max_j\|Z_j\|_{\mathrm F}\le\delta,
\]
one has
\begin{equation}
 \left(\E\mathcal H\Bigl(\sum_jU_jZ_j\Bigr)^q\right)^{1/q}
 \le\kappa.
 \label{hard:moment-bound}
\end{equation}
For $q<0$, an infinite raw moment is interpreted as a zero power mean.
The parameter $\delta$ is independent of the number of summands.
\end{lemma}

\begin{proof}
We argue by contradiction. If the assertion failed, there would be
a sequence of families violating \eqref{hard:moment-bound} whose
largest coefficient Frobenius norm tends to zero. The matrices
$\Sigma\coloneqq \sum_jZ_j^*Z_j$ are PSD of trace one, so a subsequence
converges to some $\Sigma_\infty\succeq0$ with trace one.
Each vectorized summand is centered and has deterministic norm at
most $\max_j\|Z_j\|_{\mathrm F}$. Thus the Lindeberg condition holds.
The fixed-dimensional Lindeberg central limit theorem
\cite{Billingsley1995} and
\eqref{hard:isotropy} give, along this subsequence,
\begin{equation}
 \sum_jU_jZ_j\ \Longrightarrow\
 T_\infty\coloneqq k^{-1/2}G\Sigma_\infty^{1/2},
 \label{hard:gaussian-limit}
\end{equation}
where $G$ has independent standard real or proper complex Gaussian
entries, normalized by $\E|G_{ab}|^2=1$. The row covariances are
$\Sigma_\infty/k$, distinct rows are uncorrelated, and the complex
pseudo-covariance vanishes. These assertions also allow singular
$\Sigma_\infty$.

The covariance dependence factors:
\begin{equation}
 \mathcal H(T_\infty)^2
 =|\det G|^{2/k}\det(\Sigma_\infty)^{1/k}
 \le k^{-1}|\det G|^{2/k}.
 \label{hard:covariance-factor}
\end{equation}
The inequality follows from $\tr\Sigma_\infty=1$.
Bartlett's decomposition
\cite{Muirhead1982,Goodman1963Determinant} gives independent variables
$\gamma_j\sim\operatorname{Gamma}(j\beta,\text{rate }1)$ such that
$|\det G|^2$ has the distribution of
$\beta^{-k}\prod_{j=1}^k\gamma_j$. Hence
\begin{equation}
 \E|\det G|^{2u}
 =\beta^{-ku}\prod_{j=1}^k
       \frac{\Gamma(j\beta+u)}{\Gamma(j\beta)}
 \quad(u>-\beta).
 \label{hard:bartlett}
\end{equation}
In the stated moment range,
\[
 \left[\E\left(k^{-1/2}|\det G|^{1/k}\right)^q\right]^{1/q}
 =\sqrt{\factor{k}{p}}.
\]

If $0<q<2$, \eqref{hard:quadratic-control} implies the uniform tail
bound
\[
 \E\!\left[\mathcal H(S)^q
              \mathbf1_{\{\mathcal H(S)>R\}}\right]\le R^{q-2}.
\]
Thus \eqref{hard:gaussian-limit} entails convergence of the $q$th
moments. Equation~\eqref{hard:covariance-factor} bounds the limiting
moment by $\factor{k}{p}^{q/2}<\kappa^q$.

For $q<0$, use the bounded continuous function
$\min\{\mathcal H(T)^q,M\}$, assigning it the value $M$ when
$\det T=0$. Weak convergence followed by $M\to\infty$ yields
\[
 \liminf \E\mathcal H(S)^q\ge \E\mathcal H(T_\infty)^q.
\]
If $\Sigma_\infty$ is singular, the right side is infinite.
Otherwise \eqref{hard:covariance-factor} and
\eqref{hard:bartlett} bound it below by $\factor{k}{p}^{q/2}>\kappa^q$
in the stated moment range. Taking a negative root reverses the
inequality in the required direction.

Both cases contradict the strict limiting inequalities above.
The contradiction also proves uniformity in the number of summands.
\qed
\end{proof}

\subsection{Hardness for powers below one}
\label{app:hard:subcritical}

We prove Theorem~\ref{hard:main}. Fix its field, selected dimension,
rational power $-k\beta<p<1$, $p\ne0$, and slack
$\varepsilon$, and put $q\coloneqq 2p$. The remaining powers will follow
by monotonicity. Throughout this subsection the source block length
is $b_{\rm src}$.

\subsubsection{Construction and completeness}

Write $\Pi\coloneqq \Pi_{\mathcal W}\otimes I_k$ and regard an ambient vector as blocks
$x_{i,j}\in\F^k$, indexed by $1\le i\le N$ and $1\le j\le b_{\rm src}$.
For $i$ and $\mathbf U\coloneqq (U_1,\ldots,U_{b_{\rm src}})\in\mathcal U_k^{b_{\rm src}}$, define
\begin{equation}
 L_{i,\mathbf U}x\coloneqq \sum_{j=1}^{b_{\rm src}} U_j(\Pi x)_{i,j},
 \qquad A_{i,\mathbf U}\coloneqq N L_{i,\mathbf U}^*L_{i,\mathbf U}.
 \label{hard:measurements}
\end{equation}
We call $L_{i,\mathbf U}$ the measurement map and
$A_{i,\mathbf U}$ the resulting observation. The observations are
$m\coloneqq N|\mathcal U_k|^{b_{\rm src}}$ rational PSD matrices of rank at most
$k$, acting on $\F^{N b_{\rm src}k}$. Their size is polynomial in the source
input because $b_{\rm src}$ and $k$ are fixed. Averaging first over the independent
$U_j$ and then over $i$ gives
\begin{equation}
 \frac1m\sum_{i,\mathbf U}A_{i,\mathbf U}=\Pi.
 \label{hard:average}
\end{equation}
Consequently, for every frame $P$,
\begin{equation}
 \frac1m\sum_{i,\mathbf U}\det(P^*A_{i,\mathbf U}P)^{1/k}
 \le \frac1k\tr(P^*\Pi P)\le1.
 \label{hard:one-bound}
\end{equation}
The first inequality is determinant AM--GM. Every power mean with
$p\le1$ is therefore at most one.

For a YES witness $v$ in Proposition~\ref{hard:block-source}, take
$P\coloneqq (v/\sqrt N)\otimes I_k$. This is a frame. In each block, the
coordinate selected by $v_i$ gives
$L_{i,\mathbf U}P=U_j/\sqrt N$ for some $j$, so every compression
in the YES instance is $I_k$.

\subsubsection{Soundness and the order of parameters}
\label{hard:soundness}

\emph{Reduction to source-supported matrices.}
Fix a frame $P$. Its volume scores depend only on $\Pi P$. If
$\Pi P=0$, every volume is zero. Otherwise define
\[
 Q\coloneqq \frac{\sqrt{k}\,\Pi P}{\|\Pi P\|_{\mathrm F}}.
\]
Since $\|\Pi P\|_{\mathrm F}\le\sqrt{k}$, this replacement
weakly increases every volume. Although $Q$ need not be a frame,
it is supported on $\mathcal W_{\F}\otimes\F^k$ and has Frobenius norm
$\sqrt{k}$. It is therefore enough to establish soundness on this
larger class.

Set $Z\coloneqq \sqrt{N/k}\,Q$, so that $\|Z\|_{\mathrm F}^2=N$ matches
the normalization in Lemma~\ref{tensor:hilbert}. Denote its
$k\times k$ coordinate blocks by $Z_{i,j}$. The score of $Q$ for
observation $A_{i,\mathbf U}$ is
\begin{equation}
 \det(Q^*A_{i,\mathbf U}Q)^{1/k}
 =\mathcal H\Bigl(\sum_jU_jZ_{i,j}\Bigr)^2.
 \label{hard:normalized-block-volume}
\end{equation}
Writing $S_i\coloneqq \sum_jU_jZ_{i,j}$, we can bound the $q$-mean of
$\mathcal H$ at $q=2p$ and square the result. Indeed, the same
finite average over observations $(i,\mathbf U)$ gives
\[
 \Mean_p[\mathcal H(S_i)^2]
 =\bigl(\Mean_{2p}[\mathcal H(S_i)]\bigr)^2.
\]
\emph{A common good set.}
Arrange the entries of each $Z_{i,j}$ into one row of a
coefficient matrix, with row index $(i,j)$. Its columns belong to
$\mathcal W_{\F}$ and its squared Frobenius norm is $N$.
Lemma~\ref{tensor:hilbert} supplies a common set $\mathcal G$ with
\begin{equation}
 g\coloneqq \frac{|\mathcal G|}{N}\ge1-f,\qquad
 f\coloneqq 8\tau_0+M^{-2},\qquad
 \max_j\|Z_{i,j}\|_{\mathrm F}<d_0\quad(i\in\mathcal G).
 \label{hard:good-blocks}
\end{equation}
The set may depend on $P$, but the source bound holds for every
coefficient matrix.

\emph{Moments within each block.}
Put $s_i^2\coloneqq \sum_j\|Z_{i,j}\|_{\mathrm F}^2$, so that
$\sum_i s_i^2=N$. Applying Lemma~\ref{hard:moment-lemma} requires
dividing the coefficients by $s_i$, which may be small. Introduce
a cutoff $\lambda>0$ and choose $d_0\le\delta \lambda$, where $\delta$ is
the lemma's threshold. For $s_i>\lambda$, a good block then has normalized
coefficients of norm less than $d_0/s_i\le\delta$; smaller blocks
are controlled directly by their second moment.

Expectations below average the finite gadget within each block.
On any block, homogeneity rescales \eqref{hard:quadratic-control}
from unit coefficient energy to $s_i^2$. Thus the
$q$-mean of $\mathcal H(S_i)$ is at most $s_i$ by power-mean
monotonicity (and is zero if $s_i=0$). On a
good block with $s_i>\lambda$, divide its coefficients by $s_i$ and
apply Lemma~\ref{hard:moment-lemma}: the $q$-mean is at most
$\kappa s_i$. Good blocks with $s_i\le \lambda$ have $q$-mean at most
$\lambda$.

\emph{Averaging at positive powers.}
For $0<q<2$, the good-block contribution to the raw moment is at
most $\kappa^q+\lambda^q$. Since the discarded fraction is at most $f$,
H\"older's inequality gives
\[
 \frac1N\sum_{i\notin\mathcal G}\E\mathcal H(S_i)^q
 \le\frac1N\sum_{i\notin\mathcal G}s_i^q
 \le f^{1-q/2}.
\]
Consequently
\begin{equation}
 \left(\frac1N\sum_i\E\mathcal H(S_i)^q\right)^{1/q}
 \le
 \left[\kappa^q+\lambda^q+
          f^{1-q/2}\right]^{1/q}.
 \label{hard:positive-soundness}
\end{equation}

\emph{Averaging at negative powers.}
For $q<0$, an upper bound on the power mean corresponds to a lower
bound on the raw moment. We retain only good blocks in that sum;
taking the negative root will reverse the inequality. An infinite
raw moment gives power mean zero. Put
$t_i\coloneqq \max\{\lambda,\kappa s_i\}$ on $\mathcal G$.
The block estimates now imply
$\E\mathcal H(S_i)^q\ge t_i^q$.
Discarding the other nonnegative raw moments and using the convex
decreasing function $x\mapsto x^{q/2}$ yields
\begin{align}
 \frac1N\sum_i\E\mathcal H(S_i)^q
 &\ge \frac1N\sum_{i\in\mathcal G}t_i^q \notag\\
 &\ge
 g\left(\frac{\sum_{i\in\mathcal G}t_i^2}{gN}\right)^{q/2}
 \ge g\left(\frac{\lambda^2+\kappa^2}{g}\right)^{q/2}.
 \label{hard:negative-raw}
\end{align}
After taking the negative root,
\begin{equation}
 \left(\frac1N\sum_i\E\mathcal H(S_i)^q\right)^{1/q}
 \le \sqrt{\lambda^2+\kappa^2}\,g^{1/q-1/2}.
 \label{hard:negative-soundness}
\end{equation}
Infinite moments only strengthen \eqref{hard:negative-raw}.

\emph{Choosing parameters.}
In view of \eqref{hard:normalized-block-volume}, the target
$p$-mean is the square of the $q$-mean in
\eqref{hard:positive-soundness} or \eqref{hard:negative-soundness}.
Use the following order, with analytic slack $\varepsilon/2$;
the remaining half is available for the later positive-definite
perturbation.
\begin{enumerate}
\item Choose $\kappa>\sqrt{\factor{k}{p}}$ with
$\kappa^2<\factor{k}{p}+\varepsilon/2$. Then choose $\lambda>0$ and
$0<f_*<1$ small enough that the corresponding squared bound is
below $\factor{k}{p}+\varepsilon/2$ whenever $f\le f_*$ and
$g\ge1-f_*$.
\item Obtain $\delta$ from Lemma~\ref{hard:moment-lemma}.
Choose $M\ge1$ large and $\tau_0>0$ small enough that
$8\tau_0+M^{-2}\le f_*$, and choose the coefficient threshold
$0<d_0<\min\{1,\delta \lambda\}$.
\item Invoke the source theorem with
$\tau\le\min\{\tau_0,1/2\}$, source coordinate threshold
$\delta_0\le d_0/2$, and $\eta\le(4M)^{-1}$, capping the latter
parameters at $1/2$. This fixes $b_{\rm src}$; only then does $N$ grow.
\end{enumerate}
The comparison threshold is independent of the number of coefficients,
so choosing it before the source block length creates no circular
dependence.
This proves PSD hardness strictly above the negative boundary at
nonzero power. Section~\ref{hard:pd} supplies the positive-definite
version, and Section~\ref{hard:endpoints} treats the remaining powers.

\emph{Why the arithmetic endpoint needs a new argument.}
At $q=2$, the exponent $1-q/2$ in
\eqref{hard:positive-soundness} vanishes. To illustrate this
limitation of the energy estimate, consider an exceptional fraction
$f>0$ of blocks with $s_i=f^{-1/2}$ there and $s_i=0$ elsewhere.
Then
\[
 \frac1N\sum_i s_i^2=1,\qquad
 \frac1N\sum_{i\text{ exceptional}}s_i^q=f^{1-q/2}.
\]
For $0<q<2$ this tends to zero with $f$, but for $q=2$ it equals
one. This energy profile illustrates the obstruction; it is not a
constructed source instance. The tensor construction below controls
concentrated energy by averaging several orientations while preserving
perfect completeness.

\subsection{Hardness at power one}
\label{app:hard:arithmetic}\label{tensor:appendix}

We prove Theorem~\ref{tensor:main}. Fix $k\ge2$ and abbreviate
\begin{equation}\label{tensor:constant}
 a_k\coloneqq \factor{k}{1}=\frac1k\E|\det G|^{2/k},
\end{equation}
where $G$ is a square standard Gaussian matrix.
For any fixed tensor order $d\ge2$ the construction below has
soundness at most $\max\{a_k,1/d\}$ plus an arbitrarily small
fixed slack. Choose $1/d\le a_k$ to obtain the theorem.
The tensor order is fixed before the source input grows.

An \emph{orientation} is the choice of a factor $h$ whose block index
is shared by all output rows; the other factors are sampled independently
in each row. We will use the same good sets $\mathcal G_u$ for every
orientation. After controlling exceptional independent indices, the
retained coefficients have either no bad factor or only the shared
factor bad. Write $e_\varnothing$ and $e_h$ for the corresponding
fractions of coefficient energy. The estimates below will give
\[
 F_h(T)\le a_ke_\varnothing+e_h+E_{\rm err},\qquad
 F(T)\le\max\{a_k,1/d\}+E_{\rm err}.
\]
Here $F_h$ is the arithmetic objective in orientation $h$, $F$ is
their average, and $E_{\rm err}$ can be made arbitrarily small.
Energy with exactly one bad factor is charged at full strength in
only one orientation, so averaging gives it coefficient $1/d$.
The following subsections construct the instance and justify this
accounting.

We use $u$ for a tensor factor, $h$ for the shared factor, and $r$
for an output row. Within factor $u$, $i_u$ is a source-block index
and $a_u$ a coordinate inside that block. As before, $N$ grows,
whereas the source block length $b_{\rm src}$ is fixed by the
analytic parameters.

\subsubsection{Scalar comparison on prescribed subsets}

The first ingredient replaces one independently sampled factor by
a Gaussian contraction, meaning a sum against Gaussian coefficients
along that factor. The comparison must remain valid when a prescribed
set of source blocks is retained.

Put $\theta\coloneqq 2/k\in(0,2)$ and $\mu_\theta\coloneqq \E|g|^\theta$, where $g$ is a standard
real or proper complex Gaussian. Let $\sigma_a$ be independent
uniform signs or fourth-root phases, respectively.

\begin{lemma}
\label{tensor:subset}
For every $\epsilon_S>0$, source parameters can be chosen so that,
for every $c\in \mathcal W_{\F}$ and every subset $S\subseteq[N]$,
\begin{align}
 \frac1N\sum_{i\in S}\E\left|\sum_a\sigma_ac_{i,a}\right|^\theta
 &\le
 \mu_\theta\left(\frac1N\sum_{i\in S,a}|c_{i,a}|^2\right)^{\theta/2}
 \notag\\
 &\quad+
 \epsilon_S\left(\frac1N\sum_{i,a}|c_{i,a}|^2\right)^{\theta/2}.
 \label{tensor:subset-bound}
\end{align}
The parameter choice is independent of $N$, the number of
coefficients, $c$, and $S$.
Only the full vector $c$ must belong to $\mathcal W_{\F}$; its
restriction to $S$ need not. All subset sums retain the denominator
$N$: they are averages with an indicator of $S$, not conditional
averages over $S$.
\end{lemma}

\begin{proof}
The zero vector is immediate. Otherwise normalize by homogeneity so
that $\|c\|^2=N$, and put $u_i^2\coloneqq \sum_a|c_{i,a}|^2$.
Lemma~\ref{tensor:hilbert}, applied to this one-column matrix,
gives a good set with discarded fraction at most
$f_S\coloneqq 8\tau_S+M_S^{-2}$ and coefficients bounded by $d_S$.
For a small cutoff $\lambda>0$, blocks with $u_i\le \lambda$ contribute
at most $\lambda^\theta$ to the averaged moment. The phase second-moment
identity gives $\E|\sum_a\sigma_ac_{i,a}|^\theta\le u_i^\theta$,
so H\"older bounds the discarded contribution by $f_S^{1-\theta/2}$.

On each remaining block normalize the coefficient vector by $u_i$.
The dimension-one instance of Lemma~\ref{hard:moment-lemma}, at
power $\theta/2$, uses exactly these signs or fourth-root phases.
For any $\rho>0$ it gives
\[
 \E\left|\sum_a\sigma_ac_{i,a}\right|^\theta
 \le(\mu_\theta+\rho)u_i^\theta
\]
if $d_S/\lambda$ is sufficiently small. No second central-limit
argument is required.
Restrict every bound to the prescribed set $S$ and use
\[
 \frac1N\sum_{i\in S}u_i^\theta
 \le\left(\frac1N\sum_{i\in S}u_i^2\right)^{\theta/2}.
\]
The error beyond the asserted restricted-energy term is at most
$\rho+\lambda^\theta+f_S^{1-\theta/2}$, since the total average energy is one.
Choose $\rho,\lambda,f_S$ first, then $M_S,\tau_S$ and $d_S$,
and finally the source parameters as in Lemma~\ref{tensor:hilbert}.
Rescaling gives the stated result. \qed
\end{proof}

\subsubsection{The tensor instance}

Set
\[
 \Pi\coloneqq \Pi_{\mathcal W}^{\otimes d}\otimes I_k,\qquad
 \mathcal V\coloneqq \mathcal W_{\F}^{\otimes d}\otimes\F^k.
\]
The ambient dimension is $k(Nb_{\rm src})^d$. Its coordinates have indices
\[
 ((i_1,a_1),\ldots,(i_d,a_d),t),\qquad
 i_u\in[N],\quad a_u\in[b_{\rm src}],\quad t\in[k].
\]
For a first example, take $d=2$ and share factor $h=1$. The
measurement map is
\[
 (Lx)_r\coloneqq \sum_{a,b}\sigma_b^{2,r}\,e_r^TU_a
                    (\Pi x)_{j,a,i_{2,r},b}.
\]
The first factor's block index $j$ and matrices $U_a$ are common
to all rows. The second factor's block index $i_{2,r}$ and phases
$\sigma_b^{2,r}$ are independent for each row. Orientation $h=2$
interchanges these roles. This example explains the construction;
the required tensor order is still chosen to satisfy $1/d\le a_k$.

\begin{samepage}
In general, the instance has $d$ equally weighted orientations.
For orientation $h$, use the following independent uniform choices:
\begin{enumerate}
\item Choose the shared block index $j\coloneqq i_h\in[N]$ and one
$U_a\in\mathcal U_k$ for each shared label $a\in[b_{\rm src}]$.
\item For each nonshared factor $u\ne h$ and row $r\in[k]$, choose
a block index $i_{u,r}\in[N]$ and a family of independent scalar
signs or fourth-root phases $(\sigma^{u,r}_a)_a$.
\end{enumerate}
\end{samepage}
All finite choices are enumerated in the reduction; probability
notation describes their uniform average. Define the $r$th row by
\begin{equation}
 (Lx)_r\coloneqq 
 e_r^T\sum_{a_1,\ldots,a_d}
 \left(\prod_{u\ne h}\sigma^{u,r}_{a_u}\right)U_{a_h}
 (\Pi x)_{i_{1,r},a_1,\ldots,j,a_h,\ldots,i_{d,r},a_d},
 \qquad A\coloneqq N^dL^*L.
 \label{tensor:measurement}
\end{equation}
The number of observations is a constant times $dN^{1+k(d-1)}$; hence the output
size is polynomial for fixed $k,d,b_{\rm src}$. Its entries are rational,
and each observation has rank at most $k$. Gaussian variables enter
only the analysis; no Gaussian samples are encoded in the instance.

Averaging the scalar phases kills unequal labels in nonshared
factors. Equation~\eqref{hard:isotropy} kills unequal shared labels
and internal coordinates and assigns $1/k$ to each output row.
Averaging the source indices then gives
\begin{equation}
 \E A=\Pi.
 \label{tensor:average}
\end{equation}
Thus determinant AM--GM bounds every feasible arithmetic objective
by one. For a YES witness $v$, the frame
$P\coloneqq (v^{\otimes d}/N^{d/2})\otimes I_k$ has every scaled output
$N^{d/2}LP$ equal to a diagonal phase matrix times the shared
unitary selected by $v_j$. Every compression is therefore $I_k$.

\subsubsection{Common good sets and independent-index replacement}

\emph{Good sets and coefficient energy.}
For soundness, project a candidate frame onto $\mathcal V$ and,
unless it is zero, scale the projection to a matrix $T$ with
$\|T\|_{\mathrm F}^2=k$. This can only increase the objective.
For each tensor factor $u$, reshape $\sqrt{N/k}\,T$ with row index
$(i_u,a_u)$ and all other indices as columns.
The reshaped matrix has columns in $\mathcal W_{\F}$ and squared
Frobenius norm $N$. Apply Lemma~\ref{tensor:hilbert} with
$\tau_0\coloneqq\tau_H$, where $\tau_H$ will be chosen in the
parameter list below. It supplies a set
$\mathcal G_u\subseteq[N]$ with discarded fraction at most
$f\coloneqq 8\tau_H+M^{-2}$ such that
\begin{equation}
 \sum_{a_u,\ {\rm other\ indices}}|T|^2\le\frac{kM^2}{N},
 \qquad
 \max_{a_u}\sum_{{\rm other\ indices}}|T|^2
 \le\frac{kd_0^2}{N}
 \quad(i_u\in\mathcal G_u).
 \label{tensor:good-slices}
\end{equation}
The sets are fixed by $T$, before any random choices in the measurement map.
They are the same in every orientation.

The good/bad choices in the $d$ factors partition the coefficient
tensor into disjoint rectangles, each a Cartesian product of good
or bad block-index sets. Let $e_\varnothing$ be $1/k$ times its squared
Frobenius energy where every $i_u$ is good. Let $e_h$ be the same
normalized energy where only $i_h$ is bad. Then
\begin{equation}
 e_\varnothing+\sum_{h=1}^d e_h\le1.
 \label{tensor:energy-simplex}
\end{equation}
The missing energy lies in rectangles with at least two bad factors.
These are disjoint portions of coefficient energy; the determinant
itself is not decomposed additively across rectangles.

\begin{samepage}
For $d=2$, the four energies sum to one. The table describes the main retained
branches; shared-factor and discarded-event errors are bounded below.
\begin{center}
\small
\begin{tabularx}{\linewidth}{@{}lllX@{}}
\toprule
Factor 1 & Factor 2 & Energy & Treatment \\
\midrule
Good & Good & $e_\varnothing$ & Coefficient $a_k$ in either orientation. \\
Bad & Good & $e_1$ & Coefficient one only when factor 1 is shared. \\
Good & Bad & $e_2$ & Coefficient one only when factor 2 is shared. \\
Bad & Bad & $e_{12}$ & Excluded from either main retained branch. \\
\bottomrule
\end{tabularx}
\end{center}
\end{samepage}
\emph{Discarding exceptional independent indices.}
Fix orientation $h$ and write $R\coloneqq N^{d/2}LT$ for its scaled output.
Since $A=N^dL^*L$, the observation score of $T$ is exactly
$\det(T^*AT)^{1/k}=|\det R|^{2/k}$.
Condition on the shared index and unitary family. The rows are
independent under the remaining indices and phases. Hadamard gives
$|\det R|^{2/k}\le\prod_r\|R_r\|^{2/k}$.
For a specified bad independent index in row $r$, H\"older yields
\[
 \E[\mathbf1_{\rm bad}\|R_r\|^{2/k}\mid{\rm shared}]
 \le f^{1-1/k}
       \bigl(\E[\|R_r\|^2\mid{\rm shared}]\bigr)^{1/k}.
\]
Apply the ordinary moment inequality in the other rows, then
AM--GM to their second moments. Averaging the shared data and
using $\E\|R\|_{\mathrm F}^2=k$ for the unrestricted measurement
shows that this bad event
contributes at most $f^{1-1/k}$. A union bound over
$m_0\coloneqq k(d-1)$ pairs $(u,r)$---one for each nonshared factor in
each row---shows that restricting every independent
index to $\mathcal G_u$ costs at most
\begin{equation}
 m_0 f^{1-1/k}.
 \label{tensor:bad-cost}
\end{equation}
No restriction on the energy carried by the discarded indices has
been used.

\emph{Replacing retained factors by Gaussian contractions.}
Next replace each independently indexed factor in each row by a
Gaussian contraction. Condition on all other variables. The
determinant is a scalar linear form in the current factor, and its
full coefficient vector belongs to $\mathcal W_{\F}$ because contractions
in the other tensor factors preserve that membership.
Lemma~\ref{tensor:subset} applies to this full vector while restricting
the average to $S$; it requires no closure under restriction to $S$.
For a retained set $S$, let $G_S$ contract with independent
coefficients $N^{-1/2}g_{i,a}\mathbf1_{i\in S}$, and let
$G_{\rm all}$ use all indices. The scalar Gaussian moment identity
and Lemma~\ref{tensor:subset} give
\begin{equation}
 \E_{I,\sigma}
 [\mathbf1_{I\in S}|\det R|^{2/k}]
 \le
 \E_{G_S}|\det R|^{2/k}
 +\xi\,\E_{G_{\rm all}}|\det R|^{2/k},
 \qquad \xi\coloneqq \epsilon_S/\mu_{2/k}.
 \label{tensor:one-replacement}
\end{equation}
We use this comparison conditionally for the determinant.
Every subsequent replacement has the same scalar linear form.

Iterating \eqref{tensor:one-replacement} produces a positive sum of
Gaussian branches. For example, two replacements give four terms
with coefficients $1,\xi,\xi,\xi^2$. After $m_0$ replacements the
coefficients sum to $(1+\xi)^{m_0}$. The branch retaining every
$\mathcal G_u$ has coefficient one, so the other coefficients sum
to $(1+\xi)^{m_0}-1$.
Every branch has expected volume at most one. Indeed, after all
replacements its rows have the form
\begin{equation}
 R_r=\sqrt N\,e_r^T\sum_aU_aY_{r,a},
 \label{tensor:row-chaos}
\end{equation}
where $Y_{r,a}$ contracts the slice of $T$ with $d-1$ independent
standard Gaussian vectors over retained or unrestricted sets.
The factor $\sqrt N$ results from cancelling $d-1$ covariance
factors $1/N$ against the original scale $N^d$.
If $T_r^{\rm ret}$ denotes the tensor with the coordinate restrictions
used in row $r$, phase isotropy and Gaussian second moments give
\[
 \E\|R\|_{\mathrm F}^2
 =\frac1k\sum_{r=1}^k\|T_r^{\rm ret}\|_{\mathrm F}^2
 \le k.
\]
A restriction on the shared index only deletes nonnegative terms.
This is an inequality because coordinates have now been removed;
the earlier second-moment equality concerned the unrestricted tensor.
Determinant AM--GM bounds each branch by one, proving the total
replacement error
\begin{equation}
 (1+\xi)^{m_0}-1.
 \label{tensor:expansion-cost}
\end{equation}
In the main branch every nonshared factor is restricted to its
common good set. Therefore the collections $(Y_{r,a})_a$ are
independent and identically distributed across rows.

\subsubsection{The shared factor and the surviving rectangles}

The independent-factor replacements have made each row a
multilinear function of Gaussian vectors. Positive Gaussian degree
allows us to compare the remaining shared factor by covariance
averaging, without another central-limit argument.

\begin{lemma}[Shared-factor comparison]
\label{tensor:shared}
Fix $k\ge2$, $M_1\ge1$, and $\epsilon_{\rm sh}>0$.
There is $\delta_1>0$, depending only on these parameters, with
the following property. Each $k\times k$ random matrix $X_{r,a}$
is a multilinear contraction of a deterministic coefficient tensor
with $D\ge1$ independent standard Gaussian vectors. Within a row
$r$, the same $D$ vectors are used for every label $a$.
Different rows use independent copies, and the collections
$(X_{r,a})_a$ are identically distributed across rows. Suppose
\[
 \max_a\E\|X_{r,a}\|_{\mathrm F}^2\le\delta_1^2/N,
 \qquad
 \sum_a\E\|X_{r,a}\|_{\mathrm F}^2\le M_1^2/N.
\]
For independent $U_a\sim\mathcal U_k$, independent also of the
Gaussian vectors, define
\[
 Y_{r,:}\coloneqq \sqrt N\sum_a e_r^TU_aX_{r,a},\qquad
 \Sigma\coloneqq \sum_a\E(X_{r,a}^*X_{r,a}).
\]
Then
\begin{equation}
 \E|\det Y|^{2/k}
 \le Na_k\det(\Sigma)^{1/k}+\epsilon_{\rm sh}
 \le\frac{Na_k}{k}\tr\Sigma+\epsilon_{\rm sh}.
 \label{tensor:shared-bound}
\end{equation}
The choice of $\delta_1$ is independent of $N$, the number of labels,
all Gaussian-vector dimensions, and the positive degree $D$.
\end{lemma}
\begin{proof}
\emph{Step 1: average the covariances of independent Gaussian mixtures.}
For PSD row covariances define
\[
 \Psi(S_1,\ldots,S_k)\coloneqq 
 \E\left|\det
 \begin{bmatrix}g_1S_1^{1/2}\\ \vdots\\g_kS_k^{1/2}\end{bmatrix}
 \right|^{2/k},
\]
where the $g_r$ are independent standard Gaussian rows.
If independent rows $R_r$ are mixtures of centered Gaussian rows,
with unconditional covariances $S_r\coloneqq \E(R_r^*R_r)$, then
\begin{equation}
 \E|\det R|^{2/k}\le\Psi(S_1,\ldots,S_k).
 \label{tensor:gaussian-mixture}
\end{equation}
Indeed, condition on all rows but one and expand the determinant
as a linear form in that row. If its covariance in the mixture is
$K$, the conditional moment is $\mu_{2/k}(z^*Kz)^{1/k}$ for a
fixed coefficient vector $z$. Concavity gives
$\E_K(z^*Kz)^{1/k}\le(z^*S_rz)^{1/k}$, replacing the row by
an independent Gaussian of covariance $S_r$. Repeat over the rows.
This uses row independence and separate covariance concavity,
not joint concavity of $\Psi$. The same conditioning shows that
$\Psi$ is increasing in each covariance in the PSD order.

\emph{Step 2: condition on the cyclic shifts.}
Write $U_a=\diag(\varepsilon_{a,1},\ldots,\varepsilon_{a,k})C_{t_a}$.
Condition on the cyclic shifts $\mathbf t\coloneqq (t_a)_a$, leaving the
row phases random. The rows of $Y$ are then independent, because
both their phases and their Gaussian groups are independent across
$r$. Each row is a mixture of centered Gaussian rows: condition
further on its phases and all but one Gaussian group. The remaining
dependence is linear in a standard Gaussian vector, and is proper
in the complex case. This is where $D\ge1$ is essential.

\emph{Step 3: control covariance fluctuations.}
For deterministic PSD matrices
\[
 B_{a,j}\coloneqq \E(X_{r,a}^*e_je_j^*X_{r,a}),\qquad
 b_a\coloneqq \sum_j\tr B_{a,j}=\E\|X_{r,a}\|_{\mathrm F}^2,
\]
which do not depend on $r$, phase averaging gives the conditional
row covariances
\[
 S_r(\mathbf t)=N\sum_a B_{a,\pi_{t_a}(r)},
 \qquad S_0\coloneqq \E_{\mathbf t}S_r(\mathbf t)=N\Sigma/k.
\]
By \eqref{tensor:gaussian-mixture},
$\E[|\det Y|^{2/k}\mid\mathbf t]\le\Psi(S_1(\mathbf t),\ldots,S_k(\mathbf t))$.
No individual label carries much energy. Quantitatively,
independence of the shifts across labels and
$\|B\|_{\mathrm F}\le\tr B$ for PSD $B$ give
\begin{equation}
 \begin{split}
 \E_{\mathbf t}\|S_r(\mathbf t)-S_0\|_{\mathrm F}^2
 &\le\frac{N^2}{k}\sum_{a,j}\|B_{a,j}\|_{\mathrm F}^2\\
 &\le\frac{N^2}{k}\sum_a b_a^2
 \le\frac{\delta_1^2M_1^2}{k}.
 \end{split}
 \label{tensor:covariance-fluctuation}
\end{equation}
For every choice of shifts,
$\sum_r\tr S_r(\mathbf t)=N\tr\Sigma\le M_1^2$.
Define
\[
 e(\mathbf t)\coloneqq \max_r\|S_r(\mathbf t)-S_0\|_{\rm op}.
\]
Consequently,
\[
 0\le e(\mathbf t)\le M_1^2,
 \qquad \E_{\mathbf t}e(\mathbf t)^2\le\delta_1^2M_1^2.
\]
For the first bound, the norm of a difference of PSD matrices is
at most the larger of their operator norms; for the second, use
\eqref{tensor:covariance-fluctuation} and sum over $r$.

\emph{Step 4: convert covariance error into determinant error.}
Since $S_r(\mathbf t)\preceq S_0+e(\mathbf t)I$,
\[
 \Psi(S_1(\mathbf t),\ldots,S_k(\mathbf t))
 \le k a_k\det(S_0+e(\mathbf t)I)^{1/k}.
\]
All eigenvalues of $S_0$ are at most $M_1^2/k$. Telescoping their
product after adding $e\le M_1^2$ gives
\[
 0\le\det(S_0+eI)-\det S_0\le ke(2M_1^2)^{k-1}.
\]
Apply $(x+y)^{1/k}-x^{1/k}\le y^{1/k}$ and then
$\E e^{1/k}\le(\E e^2)^{1/(2k)}$. Using $a_k\le1$, which
follows from Gaussian determinant--trace AM--GM, yields an error
that tends to zero with $\delta_1$, uniformly in the number of
labels and all contraction dimensions:
\begin{equation}
 \E|\det Y|^{2/k}
 \le Na_k\det(\Sigma)^{1/k}
   +C_k M_1^{2-1/k}\delta_1^{1/k},
 \quad C_k\coloneqq k(k2^{k-1})^{1/k}.
 \label{tensor:explicit-shared}
\end{equation}
Choose $\delta_1$ to make this error at most $\epsilon_{\rm sh}$.
Finally $\det(\Sigma)^{1/k}\le\tr\Sigma/k$ proves the second
inequality. Every estimate allows singular covariances. \qed
\end{proof}

\emph{Application to the retained tensor.}
Apply the lemma to \eqref{tensor:row-chaos} at a good shared index
$j\in\mathcal G_h$, with $D=d-1$. If $m_{j,a}$ is the squared Frobenius
norm of its tensor slice after restricting every other factor to
its good set, Gaussian second moments give
$\E\|Y_{r,a}\|_{\mathrm F}^2=m_{j,a}$.
Equation~\eqref{tensor:good-slices} gives
\[
 \sum_a m_{j,a}\le kM^2/N,\qquad
 \max_a m_{j,a}\le kd_0^2/N.
\]
Thus take $M_1\coloneqq \sqrt{k}M$ and require $\sqrt{k}d_0\le\delta_1$.
Let $\Sigma_j\coloneqq \sum_a\E(Y_{r,a}^*Y_{r,a})$.
The trace of $\Sigma_j$ is exactly the retained tensor energy at
shared block $j$. Averaging \eqref{tensor:shared-bound} over the
good shared indices gives
\begin{equation}
 \frac1N\sum_{j\in\mathcal G_h}\E|\det R|^{2/k}
 \le \frac{a_k}{k}\sum_{j\in\mathcal G_h}\tr \Sigma_j+\epsilon_{\rm sh}
 =a_ke_\varnothing+\epsilon_{\rm sh}.
 \label{tensor:good-energy}
\end{equation}
For a bad shared index, no Gaussian approximation is needed.
The second-moment identity gives
$\E\|R\|_{\mathrm F}^2=N\tr \Sigma_j$ at fixed $j$, so
\begin{equation}
 \frac1N\sum_{j\notin\mathcal G_h}\E|\det R|^{2/k}
 \le\frac1k\sum_{j\notin\mathcal G_h}\tr \Sigma_j=e_h.
 \label{tensor:bad-energy}
\end{equation}
The other factors have already been restricted to their good sets.
Accordingly, the right side is the energy where \emph{only} factor
$h$ is bad.

Combining \eqref{tensor:bad-cost},
\eqref{tensor:expansion-cost}, \eqref{tensor:good-energy},
and \eqref{tensor:bad-energy}, the objective $F_h(T)$ of orientation
$h$ satisfies
\[
 F_h(T)\le a_ke_\varnothing+e_h+E_{\rm err},\qquad
 E_{\rm err}\coloneqq m_0f^{1-1/k}+(1+\xi)^{m_0}-1+\epsilon_{\rm sh}.
\]
Averaging orientations and applying \eqref{tensor:energy-simplex}
therefore yields
\begin{equation}
 F(T)\le a_ke_\varnothing+\frac1d\sum_h e_h+E_{\rm err}
 \le\max\{a_k,1/d\}+E_{\rm err}.
 \label{tensor:final}
\end{equation}

\emph{Error budget and parameter order.}
The three terms in $E_{\rm err}$ come respectively from exceptional
independent indices, Gaussian replacement branches, and the shared
factor. For target slack $\varepsilon$, make each smaller than
$\varepsilon/6$. This gives analytic error below $\varepsilon/2$
and leaves the other half for the positive-definite perturbation.
The choices have the following order:
\begin{enumerate}
\item Fix $k,d$ and $\varepsilon$. Choose $M\ge1$ large and
$\tau_H>0$ small enough that
$m_0(8\tau_H+M^{-2})^{1-1/k}<\varepsilon/6$.
\item Choose $0<\epsilon_{\rm sh}<\varepsilon/6$ and obtain
$\delta_1$ from Lemma~\ref{tensor:shared} at $M_1\coloneqq \sqrt{k}M$.
Then choose $0<d_H<1$ with $\sqrt{k}d_H\le\delta_1$.
\item Choose $\epsilon_S>0$ so that
$(1+\epsilon_S/\mu_{2/k})^{m_0}-1<\varepsilon/6$.
The proof of Lemma~\ref{tensor:subset} then supplies
$M_S,\tau_S,d_S$.
\item Invoke the source theorem once, with parameters satisfying
\[
 \begin{aligned}
 \tau&\le\min\{\tau_H,\tau_S,1/2\},\\
 \delta_0&\le\min\{d_H/2,d_S/2,1/2\},\\
 \eta&\le\min\{(4M)^{-1},(4M_S)^{-1},1/2\}.
 \end{aligned}
\]
\end{enumerate}
These parameters give the common good sets with coefficient
threshold $d_0=d_H$ and
also the prescribed-subset comparison. Only now is $b_{\rm src}$
fixed, after which the source input size grows.
All errors are uniform in the block length and contraction
dimensions; otherwise the earlier choices could depend circularly
on the block length returned by the source theorem. Thus
$E_{\rm err}<\varepsilon/2$ is achievable in \eqref{tensor:final},
proving PSD arithmetic hardness.

\subsection{Encoding and positive-definite observations}
\label{hard:encoding}

The analytic gaps are now established. We next check that the
instances can be encoded and feasible outputs evaluated in polynomial
time, then extend the gaps to positive-definite observations and the
remaining powers. All parameters other than the source input size
are fixed in the following bounds.

\subsubsection{Rational projections and feasible outputs}

\emph{Encoding the reduction.}
The source subspace is the kernel of an integer matrix whose
entries belong to $\{0,1,-1\}$. Select linearly independent rows
to form $D$. If there are no constraints, the projection is $I$;
otherwise
\begin{equation}
 \Pi_{\mathcal W}=I-D^T(DD^T)^{-1}D.
 \label{hard:rational-projector}
\end{equation}
The entries of $DD^T$ are integers of polynomial magnitude.
The adjugate formula and Hadamard's determinant bound give
polynomial bit lengths for its inverse entries. Exact Gaussian
elimination and rational matrix multiplication compute
\eqref{hard:rational-projector} in polynomial time.
All gadget matrices in \eqref{hard:measurements} consequently have
polynomial-bit rational real and imaginary parts. Enumerating
the moment-gadget parameter tuples uses a constant factor depending
only on the fixed reduction parameters.
The tensor construction has ambient dimension $k(N b_{\rm src})^d$ and
a constant times $dN^{1+k(d-1)}$ observations. Tensor powers of
$\Pi_{\mathcal W}$ multiply the entry bit length by only a fixed factor;
the subsequent rational sums and products also have polynomial
bit length. Therefore both the dense output size and the reduction
time are polynomial. This assertion fixes $k,d$, and the
approximation slack before the input size grows. No Gaussian
random numbers are encoded in either reduction.

\emph{Representing a feasible output.}
A full-column-rank rational matrix $R\in\F^{n\times k}$ specifies
the exactly feasible subspace with orthonormal frame
$P\coloneqq R(R^*R)^{-1/2}$. Its objective can be evaluated without
representing the inverse square root:
\begin{equation}
 \det(P^*A_iP)^{1/k}
 =
 \left(\frac{\det(R^*A_iR)}{\det(R^*R)}\right)^{1/k}.
 \label{hard:rational-value}
\end{equation}
\emph{Evaluating to promise-gap accuracy.}
The determinants in the ratio are rational and can be computed
exactly; in particular, zero volumes are detected exactly.
For fixed rational nonzero $p$, algebraic-root approximation
evaluates the power mean to any required fixed additive
promise-gap accuracy in polynomial time. The rational numbers
involved have polynomial bit length, so their nonzero values
are bounded above and below by powers of two with polynomial
exponents. When $p=0$, multiply the determinant ratios and take
the $km$th root, or add their logarithms; both give the same
polynomial-time evaluation with a fixed promise gap. If a
positive-weight volume is zero, the objective is zero for
$p\le0$.

\emph{Distinguishing the promise cases.}
A deterministic improvement on the hardness factor would
therefore imply $\mathrm P=\mathrm{NP}$. A randomized polynomial-time
algorithm returning a feasible rational subspace with an improved
factor with probability at least $2/3$ would imply
$\mathrm{NP}=\mathrm{RP}$. To see the one-sided error explicitly,
choose a rational threshold strictly between the reduced NO bound
and the promised value of a successful YES output. Reject invalid
or rank-deficient output bases. For a valid basis, compute a
certified lower bound on its objective with error smaller than the
remaining fixed gap, and accept only if this lower bound exceeds
the threshold. No feasible output can pass on a NO instance,
including outputs on failed random seeds. Every successful YES
output passes. This is an RP procedure for an NP-hard promise
problem; composing with its deterministic reduction gives
$\mathrm{NP}\subseteq\mathrm{RP}$, hence equality.

At $p=1$, such interval evaluation requires only rational
determinants and bisection for their $k$th roots. One does not need
to decide equality between an arbitrary sum of algebraic numbers
and a rational threshold. The stronger sufficient assumption
$\mathrm{NP}\not\subseteq\mathrm{BPP}$ also excludes
bounded-error value-estimation improvements. The algorithms and
the gap reductions themselves do not require a complexity assumption.

\subsubsection{Continuity at the nonzero powers used in the reductions}

The perturbation below first bounds the change in each volume
score. This elementary continuity estimate converts that bound
into a change in the objective.

\begin{lemma}[Equal-weight continuity]\label{hard:continuity}
Let $p\le1$, $p\ne0$, $x\in[0,H]^m$, $y\ge0$, and
$\|x-y\|_\infty\le t\le1$. Then
\begin{equation}
\label{hard:continuity-modulus}
 |\Mean_p[x]-\Mean_p[y]|\le
 \begin{cases}
 m^{-1/p}t,&p<0,\\
 (H+1)^{1-p}t^p/p,&0<p<1,\\
 t,&p=1.
 \end{cases}
\end{equation}
\end{lemma}
\begin{proof}
For positive coordinates and $p<0$, put $z_i\coloneqq x_i^p$.
Differentiation gives
\[
 \sum_i\frac{\partial\Mean_p}{\partial x_i}
 =m^{-1/p}
 \frac{\sum_i z_i^{1-1/p}}{(\sum_i z_i)^{1-1/p}}
 \le m^{-1/p}.
\]
Integrate along a segment and extend continuously to zero coordinates.
For $0<p<1$, use $|a^p-b^p|\le|a-b|^p$ and the derivative
bound for the outer $1/p$ power on $[0,(H+1)^p]$.
The arithmetic case is immediate. \qed
\end{proof}

\subsubsection{A polynomial-bit positive-definite perturbation}
\label{hard:pd}

Fix $-k\beta<p\le1$, $p\ne0$, in the ranges covered by
one of the two constructions, and let $A_1,\ldots,A_m$ be its
PSD instance.
We will use the perturbation
$\widetilde A_i\coloneqq (A_i+t_{\rm pd}I)/(1+t_{\rm pd})$.
Adding $t_{\rm pd}I$ makes each observation positive definite;
dividing by $1+t_{\rm pd}$ preserves every identity YES compression.
It remains to bound the increase in the NO objective and choose
$t_{\rm pd}$ with polynomial bit length.

In the one-index construction,
$\|L_{i,\mathbf U}\|_{\mathrm{op}}^2\le b_{\rm src}$, so
$\|A_i\|_{\mathrm{op}}\le N b_{\rm src}$. In the tensor construction,
each row of the unprojected measurement map has squared norm
$b_{\rm src}^d$: each label tuple contributes a unit row vector,
and distinct tuples occupy disjoint ambient coordinates. Right multiplication
by $\Pi$ cannot increase the operator norm. Hence
\[
 \|L\|_{\mathrm{op}}^2\le kb_{\rm src}^d,\qquad
 \|A_i\|_{\mathrm{op}}\le k(N b_{\rm src})^d.
\]
In either case choose the corresponding integer bound $B$ on
all observation norms. For a rational
$0<t_{\rm pd}\le1$, set
\begin{equation}
 \widetilde A_i\coloneqq \frac{A_i+t_{\rm pd} I}{1+t_{\rm pd}}.
 \label{hard:pd-perturbation}
\end{equation}
Every matrix is positive definite, every YES compression
remains $I_k$, and the average matrix is the contraction
$(\Pi+t_{\rm pd} I)/(1+t_{\rm pd})$. Determinant--trace AM--GM
therefore bounds every feasible power mean by one, so the YES
value is still exactly one.

For an arbitrary frame, let
$a_i\coloneqq \det(P^*A_iP)^{1/k}\in[0,B]$.
The compression eigenvalues belong to $[0,B]$, and hence
\[
 0\le
 \det(P^*A_iP+t_{\rm pd} I)-\det(P^*A_iP)
 \le kt_{\rm pd}(B+1)^{k-1}.
\]
The inequality $(x+y)^{1/k}-x^{1/k}\le y^{1/k}$ for
$x,y\ge0$ gives the uniform volume bound
\begin{equation}
 \det(P^*(A_i+t_{\rm pd} I)P)^{1/k}\le a_i+t,
 \qquad
 t\coloneqq \bigl[kt_{\rm pd}(B+1)^{k-1}\bigr]^{1/k}.
 \label{hard:volume-modulus}
\end{equation}
Here $t_{\rm pd}$ is the perturbation of the matrix, whereas $t$
bounds the induced increase in a volume score before normalization.
Dividing by $1+t_{\rm pd}$ can only decrease the upper bound.
By monotonicity, the perturbed objective is at most
$\Mean_p[a_1+t,\ldots,a_m+t]$. Lemma~\ref{hard:continuity}, with
$H=B$, bounds its increase over the original objective.

To limit the objective increase to a prescribed fixed
$\Delta>0$, it suffices to take $t\le1$ and
\[
 t\le
 \begin{cases}
 \Delta m^{1/p},&p<0,\\[2pt]
 \bigl[p\Delta/(B+1)^{1-p}\bigr]^{1/p},&0<p<1,\\[2pt]
 \Delta,&p=1,
 \end{cases}
 \qquad
 0<t_{\rm pd}\le\frac{t^k}{k(B+1)^{k-1}}.
\]
A sufficiently small dyadic $t_{\rm pd}$ has polynomial bit length:
for fixed nonzero rational $p$, the displayed bounds involve only
fixed powers of $m$, $B+1$, and the prescribed fixed slack.
Conservative dyadic choices suffice, so no exact irrational
constant needs to be encoded. Starting with half the target slack
and assigning the other half to $\Delta$ proves the positive-definite
version in this range.

\subsubsection{Zero power and the negative boundary}
\label{hard:endpoints}

The remaining cases of Theorem~\ref{hard:main} follow by
power-mean monotonicity:
\[
 p\le p'\quad\Longrightarrow\quad\Mean_p[v_i]\le\Mean_{p'}[v_i].
\]
We use an auxiliary power $p'$ for which hardness is already proved.

\emph{The geometric power.}
For $p=0$, choose a fixed rational $0<p'<1$
so small that $\factor{k}{p'}<\factor{k}{0}+\varepsilon/2$.
Apply the PSD or positive-definite $p'$-power reduction with slack
less than $\varepsilon/2$. Its identity YES compressions have
geometric mean one, whereas every NO frame satisfies
$\Mean_0[v_i]\le\Mean_{p'}[v_i]\le\factor{k}{0}+\varepsilon$.
Thus no geometric-product continuity estimate is needed for the
positive-definite construction.

\emph{The negative boundary and lower powers.}
If the target power satisfies $p\le-k\beta$, the Gamma formula
implies $\factor{k}{p'}\to0$ as $p'\downarrow-k\beta$ from above.
Indeed, its first Gamma numerator diverges while $1/p'$ remains
negative and bounded away from zero. Choose a fixed rational
$-k\beta<p'<0$ with $\factor{k}{p'}<\varepsilon/2$ and apply its
PSD or positive-definite reduction with slack less than
$\varepsilon/2$. Perfect completeness is unchanged, and in the
NO case $\Mean_p[v_i]\le\Mean_{p'}[v_i]\le\varepsilon$.
The auxiliary power is fixed before the source input grows.
This completes both versions of Theorem~\ref{hard:main};
Section~\ref{hard:pd} also completes Theorem~\ref{tensor:main}.

\subsection{Transfer to projection observations}
\label{obs:transfer-appendix}

Projective normalization separates a rank-$k$ observation's scale
from its range and transfers approximation ratios. The hardness
sources have perfect completeness: any observation of rank below
$k$ rules out an identity YES compression and permits immediate
rejection. Otherwise every source observation has rank $k$, so the
normalization below applies without rank repair.

\begin{proof}[of Lemma~\ref{obs:normalization}]
Write $A_i=U_iD_iU_i^*$ with $U_i^*U_i=I_k$ and $D_i\succ0$.
Then
\[
 \det(P^*A_iP)=\det(D_i)|\det(U_i^*P)|^2,
\]
so $v_{A_i}(P)=d_i v_{\Pi_i}(P)$, where
$d_i=\det(D_i)^{1/k}$. For $p\ne0$, substituting
\[
 S\coloneqq \sum_iw_id_i^p,\qquad \widetilde w_i\coloneqq \frac{w_id_i^p}{S}
\]
gives
\[
 \Mean_{p,w}[v_{A_i}(P)]
 =S^{1/p}\Mean_{p,\widetilde w}[v_{\Pi_i}(P)].
\]
For $p=0$, the weights stay unchanged and the multiplier is
$\prod_i d_i^{w_i}$. Each multiplier is independent of $P$, which
is why the corresponding objectives have the same approximation
ratios before weight rounding.
A rational basis $V_i$ of the range gives the rational projector
$\Pi_i=V_i(V_i^*V_i)^{-1}V_i^*$. The pseudodeterminant $e_k(A_i)$
is a positive rational number of polynomial bit length.
Its fixed rational powers admit positive rational approximations
to any prescribed multiplicative precision.

For completeness, if each unnormalized weight is approximated
within $e^{\pm\sigma}$, each normalized weight lies within
$e^{\pm2\sigma}$ of its target. Hence, for every nonnegative score
vector, the two means differ by a factor between
$e^{-2\sigma/|p|}$ and $e^{2\sigma/|p|}$. This holds also when
the mean vanishes. Approximation ratios therefore lose at most
$e^{-4\sigma/|p|}$. The assertion at zero power is exact.
\qed
\end{proof}

\begin{proof}[of Theorem~\ref{obs:transfer}]
Suppose a polynomial-time algorithm has approximation factor
$\rho>\factor{k}{p}$ on rational rank-$k$ projection observations.
Choose a source instance from Theorem~\ref{hard:main} or
\ref{tensor:main}, with fixed soundness slack $\varepsilon>0$.
Every source observation has rank at most $k$. Check these ranks
exactly. If any has rank below $k$, reject the source instance:
it cannot be a YES instance, because an identity compression
would force rank at least $k$. Otherwise all source observations
already have rank $k$ and Lemma~\ref{obs:normalization} applies.

For $p\ne0$, choose $\varepsilon$ and a fixed rational-weight
tolerance $\sigma>0$ so small that
\[
 \rho e^{-4\sigma/|p|}>\factor{k}{p}+\varepsilon.
\]
Projective normalization, with the rational weights described
above, loses at most the displayed exponential factor in an
approximation ratio. Return the projection algorithm's subspace
in the unchanged original ambient space and evaluate its original
objective. On a successful YES run this value exceeds the source
NO bound, whereas no feasible output can do so on a NO instance.
The certified gap evaluation in Section~\ref{hard:encoding}
therefore distinguishes the source promise.
At $p=0$, normalization leaves the equal source weights unchanged
and preserves the ratio exactly; omit the exponential factor.
This proves the stated inapproximability in every nonvacuous case.
\qed
\end{proof}

This transfer preserves approximation-ratio hardness. It need not
preserve identity YES compressions, equal weights away from zero,
or a projection-valued average observation.

\section{Two-view computation and Martin guarantees}
\label{proof:exactness}

The two-view purification proof is in
Section~\ref{proof:purification}. This appendix supplies the
numerical details for Corollary~\ref{obs:two-objectives} and the
Martin consequences in Section~\ref{sec:martin}.

\subsection{Exactness for increasing aggregates}
\label{proof:two-objectives}

\begin{proof}[of Corollary~\ref{obs:two-objectives}]
A frame projection $Y\coloneqq PP^*$ is feasible on the right of
\eqref{obs:two-exact}.
Conversely, \eqref{obs:direct-connection} and
Theorem~\ref{obs:purification} replace every feasible $Y$
by a frame with both scores at least as large. Monotonicity
of $G$ proves equality. Upper semicontinuity and compactness
give the attainment assertion.

For the stated power means the relaxed objective is concave:
each determinant root is concave, and $\Mean_{p,w}$ is concave
and increasing for $p\le1$. Approximate convex optimization
followed by purification gives the algorithm.
Section~\ref{obs:purification-precision} establishes its complete
precision and running-time bounds, including stable evaluation
for rational $p$ and rational-basis output.
\qed
\end{proof}

\subsection{Precision for two-view purification}
\label{obs:purification-precision}

We prove the computational assertion of
Corollary~\ref{obs:two-objectives} using its concave PSD formulation.
The case $k=n$ is immediate: return the whole ambient space.
Otherwise let $0<\varepsilon<1/2$; larger requested errors can
be handled by using a fixed smaller tolerance.
All subspaces below are ultimately encoded by rational bases.

\paragraph{A positive scale and a spectral floor.}
Let $L$ be the complete input bit length.
Each rank-$k$ rational PSD input has positive eigenvalues between
$2^{-\operatorname{poly}(n,L)}$ and
$2^{\operatorname{poly}(n,L)}$. Its positive scale
$d_i\coloneqq e_k(A_i)^{1/k}$ has bounds of the same kind.
Write
\[
 F(Y)\coloneqq \Mean_{p,w}\left[
 d_i\det(U_i^*YU_i)^{1/k}\right],
 \qquad Y\succeq0,\quad\tr Y=k.
\]
The central point $Y_0\coloneqq (k/n)I$ has relaxed scores $d_i k/n$.
Since every power mean lies between its smallest and largest
arguments, $F(Y_0)$ is a positive lower bound on the optimum
with polynomial-bit logarithm, uniformly in $p\le1$.

Put $\theta\coloneqq \varepsilon/32$ and $a\coloneqq \theta k/n$.
Mixing an optimum with $Y_0$ in weight $\theta$ loses at most
a factor $1-\theta$, by concavity and nonnegativity.
It suffices to optimize on
\[
 Y\succeq aI,\qquad \tr Y=k.
\]
Every relaxed score there lies between $d_i a$ and $d_i$:
the lower bound follows from the spectral floor, and
determinant--trace AM--GM gives the upper bound.
The domain has a relative-interior ball of polynomial-bit radius
in the trace affine space and bounded diameter.
Lemma~\ref{alg:stability} bounds the relative first derivatives
independently of the power. Weak convex optimization with the
oracles described in Appendix~\ref{app:implementation} therefore
finds a rational feasible $Y$ with
\[
 F(Y)\ge(1-2\theta)\OPT_{p,k},
\]
using a number of bits and operations polynomial in
$n,L,\log(1/\varepsilon)$. Smaller internal tolerances followed
by a final mixture with $Y_0$ reserve a positive spectral margin
and ensure exact feasibility after rational rounding.

\paragraph{Stable evaluation for input powers.}
A rational full-rank basis $R_i$ for the range of $A_i$ gives
\[
 \det(U_i^*YU_i)=
 \frac{\det(R_i^*YR_i)}{\det(R_i^*R_i)}.
\]
Thus evaluating the objective and its derivatives requires no
irrational input frame. The determinant in the denominator and
$e_k(A_i)$ are positive rationals of polynomial bit length.
Logarithms, roots, and determinant derivatives can be evaluated
to any requested accuracy with polynomial bit cost.

For clarity, this cost remains polynomial for arbitrarily
negative rational $p$. Write the positive relaxed scores as
$e^{\ell_i}$ and let $\ell_{\min}\coloneqq \min_i\ell_i$. For $p<0$,
\[
 \log F(Y)=\ell_{\min}
 +\frac1p\log\left(\sum_iw_i e^{p(\ell_i-\ell_{\min})}\right).
\]
The exponentials are at most one and their weighted sum is
at least $w_{\min}>0$. Terms below the requested absolute
tolerance can be discarded. To control a retained exponential,
evaluate its log score to the required accuracy divided by
$\max(1,|p|)$. This adds only $O(\log(2+|p|))$ operand bits.
The normalized weights
$w_ie^{p\ell_i}/\sum_jw_je^{p\ell_j}$ used in the derivative
have the same stable evaluation.
Near zero, the centered exponential and logarithm expansions
in Appendix~\ref{app:implementation} avoid cancellation; their
additional precision is polynomial in the encoding length of $p$.
At zero the log objective is simply $\sum_iw_i\ell_i$.

\paragraph{Principal vectors to certified precision.}
Principal vectors can be computed by standard singular-value
methods \cite{BjorckGolub1973}. Here is the bit-precision
justification needed for rational input.
The squared principal cosines are the eigenvalues of the
generalized Hermitian problem
\[
 R_1^*\Pi_2R_1 z=\lambda R_1^*R_1 z,\qquad
 \Pi_2\coloneqq R_2(R_2^*R_2)^{-1}R_2^*.
\]
Its characteristic polynomial has degree at most $k$ and
rational coefficients of polynomial bit length.
After clearing denominators, standard algebraic root-separation
bounds separate distinct roots, and any root different from
zero or one from those endpoints, by
$2^{-\operatorname{poly}(n,L)}$.
Square-free factorization detects multiplicities.
Certified eigenspace computation to the required accuracy thus
has polynomial bit complexity. A basis inside a repeated
eigenspace can be chosen arbitrarily.

For $0<\lambda<1$, if $u$ is a corresponding unit principal
vector in the first space, use
\[
 v\coloneqq \Pi_2u/\sqrt{\lambda},\qquad
 f\coloneqq (v-\sqrt{\lambda}\,u)/\sqrt{1-\lambda}.
\]
The divisors have polynomial-bit lower bounds.
Eigenvalues zero and one use orthogonal and common subspaces
directly. The same reasoning applies over $\mathbb C$:
the generalized problem is Hermitian and has a real rational
characteristic polynomial.
This analysis does not require simple or well-separated input
eigenvalues as an additional assumption; it uses the separation
bounds available from their rational encoding.

\paragraph{Purification and rational output.}
Apply the block construction from the main proof to the computed
$Y$. For its exact principal blocks, every block diagonal entry
is at least $a$, every positive block trace is at least $a$,
and every block trace is at most $k$. The square-root formulas
therefore require only polynomial-bit precision.
Angular clipping can be omitted: the unclipped block vectors
already dominate both determinants.
Let $W$ be the resulting exact orthonormal frame. Then
\[
 \det(U_i^*WW^*U_i)
 \ge\det(U_i^*YU_i)\ge a^k,\qquad i=1,2.
\]
Each compression eigenvalue is at most one, so its smallest
eigenvalue is also at least $a^k$.
This supplies a polynomial-bit conditioning bound for both
output log scores, since
$\log(1/a^k)=k\log(1/a)$ is polynomial in the input and
$\log(1/\varepsilon)$.

Approximate $W$ by a rational full-rank matrix $B_{\rm out}$.
If this approximation is sufficiently accurate, its projector
$B_{\rm out}(B_{\rm out}^*B_{\rm out})^{-1}B_{\rm out}^*$
is close to $WW^*$, and each score logarithm changes by at most
$\theta$. Standard determinant perturbation bounds on matrices
with eigenvalues at least $a^k$ give such a precision with
polynomially many bits.
All intermediate principal-vector and square-root errors are
included in this same output tolerance.
Logarithmic stability from Lemma~\ref{alg:stability} transfers it
to the power mean independently of $p$. The returned subspace
therefore has value at least
\[
 e^{-\theta}(1-2\theta)\OPT_{p,k}
 \ge(1-3\theta)\OPT_{p,k}
 \ge(1-\varepsilon)\OPT_{p,k}.
\]
This proves the stated deterministic running time and relative
accuracy for a rational output basis.

\subsection{The additive Martin guarantees}
\label{proof:martin}

Use the projection identity \eqref{martin:complement-identity}
to pass to dimension $s\coloneqq \min(k,n-k)$.
Let $W_i$ be the original observation frames when $k\le n/2$,
and orthogonal-complement frames otherwise.
Rational input range bases have rational nullspace bases
computable in polynomial time, so this reduction preserves
the rational input model.

The optimal loss is finite for every finite observation family:
a standard Gaussian matrix has nonsingular compression to each
observation space almost surely, and column normalization
preserves this property.
In the effective dimension, the relaxation has the
logarithmic form
\begin{equation}\label{martin:program}
 \min_{\substack{Y\succeq0\\\tr Y=s}}
 -\sum_iw_i\log\det(W_i^*YW_i).
\end{equation}
It is the geometric-mean relaxation after the substitution
$Y\coloneqq sX$. Its determinants can be evaluated from rational range
bases; no exact irrational orthonormal input frame is required.

\begin{proof}[of Theorem~\ref{martin:main}]
Apply Theorem~\ref{alg:attainment} to the effective
$s$-dimensional projection observations at power zero.
Its factor is exactly $\factor{s}{0}$.
Taking $-s$ times the logarithm gives
\eqref{martin:dimension-free}, and complementing the output back
when necessary leaves the loss unchanged.

For the finite-dimensional refinement, apply
Theorem~\ref{round:algorithm} to the same observations.
Put $b\coloneqq \directfactor{n}{s}{0}$.
Lemma~\ref{round:density} supplies a positive rational lower bound
$c\le b$ depending only on the field. Choose factor slack
\[
 \zeta\coloneqq \frac{c\min\{\varepsilon,1\}}{2s}.
\]
The returned geometric mean is at least $b-\zeta$ times its
original optimum. Since $\zeta/b\le1/2$, the additive loss is at most
\[
 -s\log(b-\zeta)
 =-s\log b-s\log(1-\zeta/b)
 \le\Delta_{s,n}+\varepsilon.
\]
Thus $1/\zeta$ is polynomial in $s$ and
$\max(1,1/\varepsilon)$, giving the stated running time.
The explicit expression for $\Delta_{s,n}$ follows by substituting
the finite-dimensional Gamma formula at zero power.

Finally, an additive guarantee $\Delta_k-\varepsilon$, uniform
over ambient dimension, would give geometric-volume ratio
\[
 \exp(-(\Delta_k-\varepsilon)/k)>\factor{k}{0}.
\]
Theorem~\ref{obs:transfer} excludes this improvement even for
equally weighted rational projection observations.
At $k=n$, every overlap has determinant of modulus one,
so the loss is zero.
\qed
\end{proof}

\subsection{The explicit two-view Martin solution}
\label{proof:martin-two}

\begin{proof}[of Proposition~\ref{martin:two}]
Apply Theorem~\ref{obs:purification} to $Y\coloneqq PP^*$.
Its principal-block construction weakly improves both
determinants of every candidate $P$, and the optional clipping
puts each block angle in $[0,\phi_j]$.
On these block frames, the logarithmic objective separates.
Each nontrivial angle maximizes
\[
 w\log\cos\theta+(1-w)\log\cos(\phi-\theta).
\]
This function is strictly concave on its finite-loss interval.
For $0<\phi<\pi/2$, its stationary equation is
$w\tan\theta=(1-w)\tan(\phi-\theta)$.
Writing $t\coloneqq \tan\theta$ and $T\coloneqq \tan\phi$ gives
$wTt^2+t-(1-w)T=0$; its positive root, rationalized,
is \eqref{martin:angles}.
At $\phi=0$ the common direction is optimal.
At $\phi=\pi/2$, maximizing
$w\log\cos\theta+(1-w)\log\sin\theta$ gives
$\cos^2\theta=w$.
For equal weights, symmetry and strict concavity give
$\theta=\phi/2$ in every block.
\qed
\end{proof}

\end{document}